\documentclass[a4 paper, 11 pt]{amsart}

\usepackage{mathtools}
\usepackage{amssymb}
\usepackage{tikz-cd}
\tikzset{shorten <>/.style={shorten >=#1,shorten <=#1}}
\usepackage[utf8]{inputenc}
\usepackage{graphicx}
\usepackage{stmaryrd}
\usepackage{bbold}
\usepackage[margin=3cm]{geometry}
\usepackage{float}
\usepackage[utf8]{inputenc}
\usepackage{amsthm}
\usepackage{mathtools}
\usepackage{amssymb}
\usepackage{amsmath}
\usepackage{blkarray}
\usepackage{enumitem}
\usepackage{color,hyperref}
\hypersetup{
	colorlinks,
	linkcolor={red!50!black},
	citecolor={blue!50!black},
	urlcolor={blue!80!black}
}

\theoremstyle{plain}
\newtheorem{thm}{Theorem}[section]
\newtheorem{lem}[thm]{Lemma}
\newtheorem{cor}[thm]{Corollary}

\newtheorem{prop}[thm]{Proposition}

\theoremstyle{definition}
\newtheorem{example}[thm]{Example}
\newtheorem{definition}[thm]{Definition}
\newtheorem{notation}[thm]{Notation}

\newtheorem{remark}[thm]{Remark}

\author{Calum Hughes and Adrian Miranda}

\thanks{This paper is based on material partly from the second named author's Master's Thesis, while supported by the MQRES stipend, and partly done by both authors while the first named author was supported by the Dame Kathleen Ollerenshaw PhD studentship and the second named author was supported by EPSRC under grant EP/V002325/2. We are both grateful for the financial support. Both authors would like to thank Nicola Gambino for his advice while preparing this paper, and the second named author would also like to thank Steve Lack and Dominic Verity for their advice during his Master's Thesis while some of this research was conducted. Many of our proofs that the properties in $\E$ imply the properties in $\CatE$ can be found in the second named author's Master's thesis \cite{miranda2022internal}, while the converses are new. We also thank members of the Australian Category Theory Seminar for helpful discussions and historical insight.}

\address{Department of Mathematics, University of Manchester, Alan Turing building, Oxford Road, Manchester M13 9PL, United Kingdom}

\title{The elementary theory of the 2-category of small categories}

\email{calum.hughes@manchester.ac.uk \\ adrian.miranda@manchester.ac.uk }

\newcommand{\A}{\mathbb{A}}
\newcommand{\B}{\mathbb{B}}

\newcommand{\C}{\mathbb C}

\newcommand{\X}{\mathbb X}
\newcommand{\Y}{\mathbb Y}
\newcommand{\Hom}{\normalfont\text{Hom}}

\newcommand{\codisc}{\mathbf{indisc}}
\newcommand{\disc}{\mathbf{disc}}

\newcommand{\Lagr}{\mathcal{L}}
\newcommand{\E}{\mathcal{E}}
\newcommand{\R}{\mathcal{R}}

\newcommand{\s}{\mathbf{Set}}
\newcommand{\Cat}{\mathbf{Cat}}
\newcommand{\CatE}{\mathbf{Cat(\E)}}

\begin{document}

\begin{abstract}
    We give an elementary description of $2$-categories $\mathbf{Cat}\left(\mathcal{E}\right)$ of internal categories, functors and natural transformations, where $\mathcal{E}$ is a category modelling Lawvere's elementary theory of the category of sets (ETCS). This extends Bourke's characterisation of $2$-categories $\mathbf{Cat}\left(\mathcal{E}\right)$ where $\mathcal{E}$ has pullbacks to take account for the extra properties in ETCS, and Lawvere's characterisation of the (one-dimensional) category of small categories to take account of the two-dimensional structure. Important two-dimensional concepts which we introduce include $2$-well-pointedness, full-subobject classifiers, and the categorified axiom of choice. Along the way, we show how generating families (resp. orthogonal factorisation systems) on $\mathcal{E}$ give rise to generating families (resp. orthogonal factorisation systems) on $\mathbf{Cat}(\mathcal{E})_{1}$, results which we believe are of independent interest. 
\end{abstract}

	\maketitle
 \tableofcontents

\section{Introduction}

  Lawvere's Elementary Theory of the Category of Sets (hereafter ETCS) \cite{lawvere1964elementary} provides a set theory which axiomatises the properties of function composition rather than those of a global set membership relation. It provides an important fragment of a category-theoretic foundation of mathematics, but is strictly weaker than the traditional foundation of mathematics given by Zermelo Fraenkel Set Theory with the Axiom of Choice (hereafter ZFC). Precisely, ZFC is equiconsistent with ETCS augmented with the axiom schema of replacement \cite{osius1974categorical}.

  In his PhD thesis \cite{lawvere1963functorial}, Lawvere also gave an elementary, first order axiomatisation of the category of categories and functors. He later advocated for the first order theory of the category of categories as a foundation of mathematics (CCAF) \cite{lawvere1966category}. In an address at the 2015 Category Theory conference in Aveiro, he called for an ``improved axiomatisation" to an explicit formulation of the principles of category theory \cite{Lawvere2016Grothendieck}. Our work is a step towards this goal, re-expressing Lawvere's foundational framework as one for category theory rather than one for set theory.
  
  In this paper, we propose a different categorification of ETCS which captures the natural two-dimensional structure of the $2$-category of small categories. This is the elementary theory of the $2$-category of small categories (ET2CSC) of the title. Our main result establishes that the theory of such $2$-categories is `Morita biequivalent' with ETCS, meaning that the two theories have biequivalent $2$-categories of models.

 ETCS lacks the expressive power needed to support certain important set theoretical constructions, such as transfinite recursion. Nonetheless, it does support many of the set theoretic constructions that most mathematicians use in everyday practise. Indeed, Lawvere's aim in giving the definition was to capture more closely those aspects of set theory which are more broadly used. It is a structuralist foundation, which prioritises the perspective of how sets relate to one another, rather than a materialist one such as ZFC which prioritises how sets are built, such as via well-founded trees. While philosophical considerations are not the focus of this paper, a reader interested in these matters should consult Chapters 1 and 5 of \cite{landry2017categories}, and the references therein. ET2CSC clarifies the position of the ordinary theory of small categories within Street's programme towards a formal category theory \cite{street1980cosmoi, street2006elementary}. It facilitates a structuralist framework in which many simple category theoretical constructions can be performed, just as ETCS does for many simple set theoretical constructions. In follow up work \cite{HughesMiranda20242categories}, we extend the present axiomatisation of the $2$-category of small categories by adding a discrete opfibration classifier that satisfies a categorified version of the axiom of replacement. This provides a $2$-dimensional analogue of categories of small maps \cite{joyal1995algebraic}, and extends the present theory to encompass ZFC and facilitate more sophisticated categorical constructions.

\subsection{Outline of main results}

Our main contribution is giving an elementary theory for the $2$-category of small categories, and showing that the $2$-category of models for this theory is biequivalent to that for Lawvere's elementary theory of the category of sets, as recalled in Definition~\ref{definition ETCS}, to follow.

\begin{definition}\label{definition ETCS}
    (\cite{lawvere1964elementary}) A category $\E$ is said to \emph{model the elementary theory of the category of sets} if the following conditions are satisfied.
\begin{enumerate}
    \item $\E$ has finite limits.
    \item $\E$ is cartesian closed.
    \item The terminal object $\mathbf{1}$ is a generator for $\E$, as recalled in Definition~\ref{definition generator} part (1).
    \item $\E$ has a natural numbers object, as recalled in Definition~\ref{definition NNO} part (1).
    \item $\E$ has a subobject classifier, as recalled in Definition~\ref{def: full subobject classifier} per the discussion in Remark~\ref{full subobjects generalise subobjects}.
    \item $\E$ satisfies the external axiom of choice, as recalled in Definition~\ref{def: External axiom of choice}.
\end{enumerate}
\end{definition}

 See \cite{leinster2014rethinking} for a gentle introduction to ETCS, and \cite{lawvere2005elementary} for technical details. For Definition~\ref{definition ETCS} part $n \in \{1, ..., 6\}$, Section $n+2$ exhibits a condition on the $2$-category $\CatE$ that is equivalent to the condition on $\E$ listed as axiom $n$ above. In particular, the main results of each of these sections are Theorem~\ref{prop finite weighted limits in Cat E}, Theorem~\ref{Theorem cartesian closedness iff}, Theorem~\ref{Theorem generator in E vs generator in CatE}, Theorem~\ref{thm NNO iff}, Theorem~\ref{Theorem subobject classifier in E iff full subobject classifier in Cat E}, and Theorem~\ref{Theorem external AC iff every ff + left orhtogonal to full mono is a left adjoint left inverse equivalence}. We collate these results in Theorem~\ref{thm:ET2CSC} to characterise up to $2$-equivalence those $2$-categories which are of the form $\CatE$ for $\E$ a model of ETCS. This is expressed in terms of the elementary theory of the $2$-category of small categories, which we introduce in Definition~\ref{def:ET2CSC}. Theorem~\ref{theorem characterisation of ET2CSC morphisms} builds upon this result to characterisation to morphisms of models, and finally Theorem \ref{Theorem morita equivalence between ETCS and ET2CSC} establishes the biequivalence between the $2$-categories of models of ETCS and ET2CSC.

\subsection{Key ideas and techniques}

\subsubsection{Internal Category Theory and Bourke's characterisation of $\CatE$}

 Section~\ref{sec:internalcategories} establishes our notation and conventions in internal category theory, and catalogues various concepts that will be used in constructions and proofs. Specifically, Subsection~\ref{subsection internal categories and functors} describes internal categories, functors and natural transformations via their truncated nerves, and also describes the $2$-category structure that these data comprise. In Subsection~\ref{subsection adjunctions between E and Cat E} we catalogue the various adjunctions between $\E$ and $\CatE_{1}$ that will be used throughout this paper.

 Sections~\ref{Section Finite Limits} (resp.~\ref{Section Cartesian Closedness}) review the well known relationships between finite limits in $\E$ and finite $2$-limits in $\CatE$ (resp. cartesian closedness of $\E$ and cartesian closedness of $\CatE$). Our work relies heavily on Bourke's characterisation up to $2$-equivalence of $2$-categories of the form $\CatE$ for $\E$ with pullbacks, recalled in Proposition~\ref{Bourke characterisation of Cat E}. We thereafter allow ourselves to assume that $\mathcal{K}$ is of this form, focusing on characterising the remaining aspects of ETCS.

\subsubsection{Generating families}

 The following new results in Section~\ref{Section Generators} are important stepping stones.

\begin{itemize}
    \item Lemma~\ref{extensive coproduct Cat(E)} shows that $\E$ has extensive coproducts if and only if $\CatE$ does.
    \item Theorem~\ref{thm:copower} part (2) shows that if in addition to the previous point $\E$ is also cartesian closed, then the $2$-category $\CatE$ also has copowers by $\mathbf{2}$.
\end{itemize}  

 As well as simplifying subsequent proofs by allowing two-dimensional aspects of limit like universal properties to be deduced from their one-dimensional counterparts, copowers by $\mathbf{2}$ are used to construct generators in $\CatE$ from those in $\E$. This is shown in Corollary~\ref{cor:generator}, a result that we think is of independent interest. Definition~\ref{def:2-well-pointed} introduces a definition of a $2$-category $\mathcal{K}$ being $2$-well-pointed. This is a two-dimensional analogue of well-pointedness for categories, and is a novel concept.

\subsubsection{Adjunctions and full subobject classifiers}

Sections~\ref{Section NNO} (resp.~\ref{Section Subobject Classifiers}) relate natural numbers objects (resp. subobject classifiers) in $\E$ to their appropriate counterparts in $\CatE$. The proofs in Sections~\ref{Section Generators},~\ref{Section NNO} and~\ref{Section Subobject Classifiers} use routine calculations involving the adjunctions $\Pi_{0}\dashv \mathbf{disc} \dashv {\left(-\right)}_{0} \dashv \mathbf{indisc}$, which are reviewed in Subsection~\ref{subsection adjunctions between E and Cat E}. Section~\ref{Section Subobject Classifiers} introduces the definition of a full subobject classifier, which is a different two-dimensional analogue of a subobject classifier to the discrete opfibration classifiers of \cite{weber2007yoneda}. These are a new concept, introduced in Definition~\ref{def: full subobject classifier}.

\subsubsection{Orthogonal factorisations and the categorified axiom of choice}

In Section~\ref{Section Axiom of Choice} we first give a condition on $2$-categories of the form $\mathcal{K}:=\CatE$ which is equivalent to the external axiom of choice in $\E$, and then re-express this condition in $2$-categorical terms without relying on being able to recognise $\mathcal{K}$ as $\CatE$. The internal formulation involves fully-faithfulness and the condition of being an epimorphism on objects. Whilst the first of these properties can be recognised representably in any $2$-category, the second cannot. Although we could appeal to Proposition~\ref{Bourke characterisation of Cat E} to content ourselves with recognising it via $\mathcal{K} \simeq \CatE$, we show that epimorphism on objects internal functors are characterised by a left orthogonality property against a representably defined class of maps $\R'$. This follows from Proposition~\ref{prop:liftfactorisation}, also of independent interest, in which we show that orthogonal factorisation systems $(\Lagr, \R)$ on $\E$ give rise to orthogonal factorisation systems $(\Lagr', \R')$ on $\CatE$. Indeed, $\R'$ is precisely the full subobjects, for which classifiers are examined in Section~\ref{Cat(E) full subobject}.

\section{Notation, conventions and background on internal category theory}\label{sec:internalcategories}

 In this section we establish the notation, terminology and conventions used in this paper, and catalogue concepts from internal category theory that will be crucial for our proofs. 

\begin{notation}
    In this paper we will use the following conventions for font.
    
    \begin{itemize}
        \item Calligraphic font $\mathcal{E}$, $\mathcal{C}$, $\mathcal{K}$ will be used for categories or $2$-categories, with the letter $\mathcal{K}$ typically being reserved for $2$-categories.
        \item Ordinary mathematical font will be used for objects in categories or in $2$-categories. These will typically be capitals $X, Y, Z$ when they are objects, and lower case $f, g, h$ when they are morphisms. Greek letters will typically be used for $2$-cells.
        \item Blackboard bold $\A, \B, \C$ will be used for internal categories. When we need to be even more careful in distinguishing data in $\CatE$ from data in $\E$, the former will be either underlined or overlined. As an example, in Definition~\ref{definition internal natural transformation} we distinguish between the $2$-cell $\overline{\alpha}: f \Rightarrow g$ in $\CatE$, and its components assigner, which is a morphism $\alpha: A_{0} \to B_{1}$ in $\E$.
    \end{itemize}
\end{notation}

\begin{remark}\label{2-categorical conventions}
    We assume some familiarity with $2$-category theory and basic notions from elementary topos theory. We briefly remind the reader of common $2$-categorical notions and conventions that are used in this paper. For general background on two-dimensional category theory, see \cite{lack20092, johnson20212}. We will refer to $\Cat$-enriched (co)limits as \emph{$2$-(co)limits}. We assume familiarity with the notion of \emph{powers} and \emph{copowers} by the category $\mathbf{2}:= \{\begin{tikzcd}
    \bullet \arrow[r] & \bullet
    \end{tikzcd}\}$, elsewhere also called cotensors and tensors by $\mathbf{2}$, respectively. A notion of \emph{finiteness} for weights for $2$-(colimits) is described in \cite{street1976limits}, and all $2$-(co)limits that we will consider are finite in this sense. We will call an adjoint equivalence in the $2$-category $\mathcal{V}$-$\mathbf{Cat}$ for $(\mathcal{V}, \otimes, I) = (\mathbf{Cat}, \times, \mathbf{1})$ a \emph{$2$-equivalence}. If $\mathcal{K}$ is a $2$-category, then $\mathbf{Disc}(\mathcal{K})$ will denote its category of discrete objects, with an object $X \in \mathcal{K}$ being called \emph{discrete} if any $2$-cell into $X$ is an identity. Note that $\mathbf{Disc}(\mathcal{K})$ is distinct from the un-capitalised $\mathbf{disc}: \E \to \CatE_{1}$, to be recalled in Remark~\ref{Disc left adjoint counit}, which sends an object to a discrete internal category. A functor (resp. $2$-functor) will be said to \emph{preserve} some structure if it does so up to isomorphism.
\end{remark}

\subsection{Internal categories and the $2$-category $\CatE$}\label{subsection internal categories and functors}

 Internal categories were formally introduced by Grothendieck in  \cite{grothendieck1960techniques}, but their structure was already implicit in \cite{ehresmann1959categories} and further early applications to differential geometry appeared in the subsequent \cite{ehresmann1963categories}. See chapter 8 of \cite{borceux1994handbook} for a modern textbook account of internal category theory, and B2 of \cite{johnstone2002sketches} for its relation to topos theory.

 Let $\Delta$ denote the skeleton of the `simplex category', whose objects are non-empty finite ordered sets and morphisms are order preserving functions. Identify each object in $\Delta$ with its representative ordered set $[n]:= \{0, 1, 2..., n\}$. For $k \leq n$, let $\delta_{k}^{n}: [n] \rightarrow [n+1]$ denote the unique monotonic function whose image does not contain $k \in [n+1]$ and let $\sigma_{k}^{n}: [n+1] \to [n]$ denote the unique monotonic function mapping two elements to $k$ and one element to every other possible output. Let $\Delta_{\leq 3}$ denote the full-subcategory of $\Delta$ on the objects $[n]$ for $0 \leq n \leq 3$.
	
\begin{definition}\label{internal category short definition}
		A category \emph{internal} to a locally small category $\mathcal{E}$ is a diagram in $\E$ as displayed below left, which sends the pushout squares in $\Delta_{\leq 3}$ displayed below right to pullback squares in $\mathcal{E}$. 
		
	$$\begin{tikzcd}
	    \Delta_{\leq 3}^\text{op} \arrow[rr, "\mathbb{C}"] && \mathcal{E}&{}
	\end{tikzcd}
			\begin{tikzcd}[column sep = 12, row sep = 12]
				n+2 \arrow[ddrr, phantom, "\lrcorner", very near start]
				&& n+1 \arrow[ll, "{\delta_{2}^{n+1}}"']
				\\
				\\
				n+1 \arrow[uu,"\delta_{0}^{n+1}"]
				&& n \arrow[ll,"\delta_{1}^{n}"] \arrow[uu,"\delta_{0}^{n}"']
			\end{tikzcd}$$
\end{definition}

\begin{remark}\label{internal category long definition}
		We unpack this definition, and establish notation and terminology which we will use in this paper. A category $ \mathbb{C} := \left(C_{0}, C_{1}, d_{0}, d_{1}, i, m\right) $ internal to $\mathcal{E}$ is given by the datum of a diagram in $\E$ as displayed below.
  
  $$\begin{tikzcd}[column sep = 12, row sep = 12]
			C_2 \arrow[rr, "m"description]\arrow[rr, shift left = 3, "\pi_{0}"]\arrow[rr, shift right = 3, "\pi_{1}"']
			&& C_1 \arrow[rr, "d_0", shift left=3] \arrow[rr, "d_1"', shift right=3]
			&& C_0\arrow[ll, "i"{description}]
		\end{tikzcd}$$ 
  
   The objects ${C}_{0}$, $C_{1} \in \E$ are called the \emph{object of objects} and \emph{object of arrows} respectively, and the morphisms $d_{1}, d_{0}, i, m$ are called \emph{source}, \emph{target}, \emph{identity assigner} and \emph{composition}. The \emph{object of composable $n$-tuples} $C_{n}$ for $n\in \{2, 3\}$ are pullbacks as depicted below.
	\begin{equation*}
		\begin{tikzcd}[column sep = 12, row sep = 12]
			{C}_{2} \arrow[ddrr, phantom, "\lrcorner", very near start] \arrow[rr, "{\pi}_{0}"] \arrow[dd, "{\pi}_{1}"']
			&& {C}_{1} \arrow[dd,"d_{1}"]\\
			\\
			{C}_{1} \arrow[rr, "d_{0}"']
			&& {C}_{0} &{}
		\end{tikzcd}\begin{tikzcd}[column sep = 15, row sep = 12]
			{C}_{3} \arrow[ddrr, phantom, "\lrcorner", very near start] \arrow[rr, "{\pi}_{3, 0}"] \arrow[dd, "{\pi}_{3, 1}"']
			&& {C}_{2} \arrow[dd,"{\pi}_{1}"]\\
			\\
			{C}_{2} \arrow[rr, "{\pi}_{0}"']
			&& {C}_{1}
   \end{tikzcd}
        \end{equation*}

  These data are subject to axioms asserting the commutativity of the diagrams displayed below. 

		\begin{itemize}
			\item Sources and targets for identities and composites:
		
		$$	 \begin{tikzcd}[column sep = 12, row sep = 12]
				{C}_{0} \arrow[rr, "i"] \arrow[rrdd,"{1}_{{C}_{0}}"']
				&& {C}_{1} \arrow[dd, "d_{0}"] &{}
				\\
				\\
				&& {C}_{0}
			\end{tikzcd} \begin{tikzcd}[column sep = 12, row sep = 12]
				{C}_{0} \arrow[rr, "i"] \arrow[rrdd,"{1}_{{C}_{0}}"']
				&& {C}_{1} \arrow[dd, "d_{1}"]\\
				\\
				&& {C}_{0} &{}
			\end{tikzcd} \begin{tikzcd}[column sep = 12, row sep = 12]
				{C}_{2} \arrow[rr, "m"] \arrow[dd, "{\pi}_{0}"']
				&& {C}_{1} \arrow[dd,"d_{0}"] & {}
				\\
				\\
				{C}_{1} \arrow[rr, "d_{0}"']
				&& {C}_{0}
			\end{tikzcd}\begin{tikzcd}[column sep = 12, row sep = 12]
				{C}_{2} \arrow[rr, "m"] \arrow[dd, "{\pi}_{1}"']
				&& {C}_{1} \arrow[dd,"d_{1}"]\\
				\\
				{C}_{1} \arrow[rr, "d_{1}"']
				&& {C}_{0}
			\end{tikzcd}$$
			
			\item The associativity and left and right unit laws for composition:
			
		$$	     \begin{tikzcd}[column sep = 12, row sep = 12]
				{C}_{3}\arrow[rrr, "m_{0}"] \arrow[dd, "m_{1}"']
				&&& {C}_{2} \arrow[dd,"m"]\\
				\\
				{C}_{2} \arrow[rrr, "m"']
				&&& {C}_{1} &{}
			\end{tikzcd}\begin{tikzcd}[column sep = 12, row sep = 12]
			{C}_{1} \arrow[rrr, "i_{0}"] \arrow[rrrdd, "1_{C_{1}}"']
			&&& {C}_{2} \arrow[dd, "m"]
			&&& {C}_{1} \arrow[lll, "i_{1}"'] \arrow[llldd,"1_{C_{1}}"]\\
			\\
			&&& {C}_{1}
		\end{tikzcd}
		$$
		Where the morphisms $m_{0}: = \left(m\pi_{3, 0}, \pi_{1}\pi_{3, 1}\right)$, $m_{1}:= \left(\pi_{0}\pi_{3, 0}, m\pi_{3, 1}\right)$,  $i_{0}:= \left(id_{0}, 1_{C_{1}}\right)$ and $i_{1}:= \left(1_{C_{0}}, id_{1}\right)$ are induced by the universal property of $C_{2}$ as a pullback. For example, the equation required for $m_{0}$ to be well-defined is witnessed by the following calculation.
  $$d_{1}.m.\pi_{3{,}0}= d_{1}.\pi_{1}.\pi_{3{,}0}= d_{1}.\pi_{0}.\pi_{3{,}1} = d_{0}.\pi_{1}.\pi_{3{,}1}$$
  \end{itemize}

These conditions correspond to the simplicial identities which must be preserved by functoriality of $\mathbb{C}: \Delta_{\leq 3}^\text{op} \to \E$.
  
\end{remark}

\begin{definition}\label{definition internal functor}
    Let $\mathcal{E}$ be a category with pullbacks and let $\mathbb{A}, \mathbb{B}: \Delta_{\leq 3}^\text{op} \rightarrow \mathcal{E}$ be categories internal to $\mathcal{E}$. An \emph{internal functor} from $\mathbb{A}$ to $\mathbb{B}$ is a natural transformation as depicted below.

    $$\begin{tikzcd}
        \Delta_{\leq 3}^\text{op}\arrow[rr, bend left = 30, "\mathbb{A}"{name=A}]\arrow[rr, bend right = 30, "\mathbb{B}"'{name=B}] && \mathcal{E}
        \arrow[from=A, to=B, Rightarrow, shorten = 5, "f"']
    \end{tikzcd}$$
\end{definition}

\begin{remark}\label{Explicit internal functor}
    Internal functors can also be defined explicitly as given by a \emph{component on objects} ${f}_{0}: A_0 \rightarrow B_0$ and a \emph{component on arrows} ${f}_{1} : {A}_{1} \rightarrow {B}_{1}$ in $\mathcal{E}$ which satisfy the commutativity of the diagrams shown in~\ref{axioms for internal functor}. Here the morphism $f_{2}:= \left(f_{1}\pi_{0}, f_{1}\pi_{1}\right)$, is induced by the universal property of $B_{2}$, as witnessed by the following calculation

    $$d_{1}.f_{1}.\pi_{0}= f_{0}.d_{1}.\pi_{0}= f_{0}.d_{0}.\pi_{1}=d_{0}.f_{1}.\pi_{1}$$

 The component $f_{3}: A_{3} \to B_{3}$ is uniquely determined from this information by the universal property of $B_{3}$ in a similar way. The diagrams below express $f$'s respect for sources, targets, identities, and composition, and they all correspond to naturality conditions for $f: \mathbb{A} \to \mathbb{B}$.

\begin{equation*}\label{axioms for internal functor}
\begin{tikzcd}[column sep = 12, row sep = 12]
				{A}_{1} \arrow[rr, "{f}_{1}"] \arrow[dd,"d_{0}^{\mathbb{A}}"']
				&& {B}_{1} \arrow[dd,"d_{0}^{\mathbb{B}}"]\\
				\\
				{A}_{0} \arrow[rr,"{f}_{0}"']
				&& {B}_{0} & {}
			\end{tikzcd} \begin{tikzcd}[column sep = 12, row sep = 12]
				{A}_{1} \arrow[rr, "{f}_{1}"] \arrow[dd,"d_{1}^{\mathbb{A}}"']
				&& {B}_{1} \arrow[dd,"d_{1}^{\mathbb{B}}"]\\
				\\
				{A}_{0} \arrow[rr,"{f}_{0}"']
				&& {B}_{0} & {}
			\end{tikzcd} \begin{tikzcd}[column sep = 12, row sep = 12]
				{A}_{0} \arrow[rr, "{f}_{0}"] \arrow[dd,"i^\mathbb{A}"']
				&& {B}_{0} \arrow[dd,"i^\mathbb{B}"] &{}\\
				\\
				{A}_{1} \arrow[rr,"{f}_{1}"']
				&& {B}_{1} &{}
			\end{tikzcd}\begin{tikzcd}[column sep = 12, row sep = 12]
				{A}_{2} \arrow[rr,"{f}_{2}"]
				\arrow[dd,"m^\mathbb{A}"']
				&& {B}_{2} \arrow[dd,"m^\mathbb{B}"]
				\\
				\\
				{A}_{1} \arrow[rr,"{f}_{1}"']
				&& {B}_{1}
\end{tikzcd}
\end{equation*}

	 The morphism $f_{2}$ is thought of as taking a composable pair in $\mathbb{A}$ and returning the composable pair given by its image under $f$. Given $\left(x, y\right): X \rightarrow A_{2}$, the morphism $f_{2}$ composes with $\left(x, y\right)$ to give $\left(f_{1} x, f_{1}y\right)$, and so the equation $f_{1}m\left(x, y\right) = m\left(f_{1}x, f_{1}y\right)$ follows by respect for composition.
\end{remark}

 \begin{remark}\label{Remark nerve}
     It is evident from their definition that internal categories and internal functors form a category, in fact a full subcategory of $[\Delta_{\leq 3}^\text{op}, \mathcal{E}]$. We write this category as $\mathbf{Cat}\left(\mathcal{E}\right)_{1}$, using the subscript `$1$' to distinguish it from the $2$-category $\CatE$ which we will recall in Proposition \ref{2 category Cat(E) theorem}. In particular, $\mathbf{Cat}\left(\mathcal{E}\right)_{1}$ is small (resp. locally small) if $\mathcal{E}$ is small (resp. locally small), since $\Delta_{\leq 3}^\text{op}$ is certainly small. The inclusion functor $N: \mathbf{Cat}\left(\mathcal{E}\right)_{1} \hookrightarrow [{\Delta}_{\leq 3}^{\text{op}}, \mathcal{E}]$, which sends an internal category to its underlying truncated simplicial object in $\mathcal{E}$, is called the \emph{nerve}.
 \end{remark}

\begin{prop}\label{arr faithful}
     Consider the functors $(-)_{0}, (-)_{1}: \CatE_{1} \to \E$, which send an internal category to its object of objects and object of arrows respectively.
		\begin{enumerate}
		    \item $(-)_{1}: \CatE_{1} \to \E$ is faithful.
                \item $(-)_{0}$ and $(-)_{1}$ preserve and jointly reflect limits.
		\end{enumerate}
\end{prop}

\begin{proof}
     For part (1), let $f, g: \mathbb{A} \rightarrow \mathbb{B}$ be internal functors in $\mathcal{E}$ such that $f_{1} = g_{1}$. We need to show that $f = g$. Since $f_{1} = g_{1}$, in particular $f_{1}i^\mathbb{A} = g_{1}i^\mathbb{A}$. Since $f$ and $g$ both preserve identities, this is equivalent to saying that $i^\mathbb{B}f_{0} = i^\mathbb{B}g_{0}$. But by sources (or targets) for identities in $\mathbb{B}$, we may compose these equal morphisms in $\mathcal{E}$ with the source (or target) map of $\mathbb{B}$ to see that $f_{0} = g_{0}$. For part (2), it is standard that the family of functors $(-)_{n}:[\Delta_{\leq 3}^\text{op}, \E] \to \E$ for $n \leq 3$ preserve and jointly reflect limits, and that limits in $\CatE_{1}$ are computed in $[\Delta_{\leq 3}^\text{op}, \E]$. But since the outputs for $n \in \{0, 1\}$ are enough to determine the rest of an internal category structure, it follows that $(-)_{0}$ and $(-)_{1}$ also jointly reflect limits.
\end{proof}
 
We now review how $\CatE_1$ can be upgraded to a $2$-category by incorporating the internal natural transformations of Definition~\ref{definition internal natural transformation}, to follow. 

  \begin{definition}\label{definition internal natural transformation}
     Given internal functors $\left(f_{0}, f_{1}\right), \left(g_{0}, g_{1}\right): \mathbb{A} \rightarrow \mathbb{B}$, an \emph{internal natural transformation}
     
     $$\begin{tikzcd}
        \mathbb{A}\arrow[rr, bend left = 30, "f"{name=A}]\arrow[rr, bend right = 30, "g"'{name=B}] && \mathbb{B}
        \arrow[from=A, to=B, Rightarrow, shorten = 5, "\overline{\alpha}"']
    \end{tikzcd}$$
     
      is a morphism $\alpha: A_{0} \rightarrow B_{1}$ called the \emph{component assigner}, making the following diagrams in $\E$ commute.
		
		\begin{itemize}
			\item Assignation of components:  the commutative diagrams displayed below left and below centre commutes.

		\item Internal naturality: the square displayed below right commutes, where the morphisms $\alpha_{0}:=\left(\alpha  d_{1}, g_{1}\right): A_{1} \rightarrow B_{2}$ and $\alpha_{1}:= \left({f}_{1}, \alpha d_{0}\right): A_{1} \rightarrow B_{2}$ are induced by the universal property of $B_2$.

  \begin{equation*}
      \begin{tikzcd}[column sep = 12, row sep = 12]
				{A}_{0} \arrow[rr, "\alpha"] \arrow[rrdd,"{f}_{0}"']
				&& {B}_{1} \arrow[dd, "d_{1}"] &{}
				\\
				\\
				&& {B}_{0}
			\end{tikzcd} 
   \begin{tikzcd}[column sep = 12, row sep = 12]
				{A}_{0} \arrow[rr, "\alpha"] \arrow[rrdd,"{g}_{0}"']
				&& {B}_{1} \arrow[dd, "{d}_{0}"]\\
				\\
				&& {B}_{0}
 \end{tikzcd}\qquad
 \begin{tikzcd}[column sep = 12, row sep = 12]
			{A}_{1} \arrow[rrr, "\alpha_{0}"] \arrow[dd,"\alpha_{1}"']
			&&& {B}_{2} \arrow[dd,"m"]\\
			\\
			{B}_{2} \arrow[rrr,"m"']
			&&& {B}_{1} & {}
\end{tikzcd}
  \end{equation*}
		\end{itemize}
 \end{definition}

Internal natural transformations correspond to simplicial homotopies $\{\alpha_{0, ..., n}: A_{n} \to B_{n+1}\}_{n \in \mathbb{N}}$ \cite{goerss2009simplicial}, but are once again determined by significantly less data than in the setting of general simplicial objects due to the universal property of pullbacks in $\E$.

\begin{prop}\label{2 category Cat(E) theorem}
 	(Proposition 8.1.4 of \cite{borceux1994handbook}, Section 1.4 of \cite{miranda2022internal}) Let $\mathcal{E}$ be a category with pullbacks. Categories, functors and natural transformations internal to $\mathcal{E}$ form a $2$-category $\mathbf{Cat}\left(\mathcal{E}\right)$ whose underlying category is $\CatE_{1}$, identity $2$-cells $\overline{1_{f}}$ have component assigners given by $if_{0}$, vertical composite of $2$-cells below left has component assigner given by the morphism in $\E$ depicted below right. 
  
$$\begin{tikzcd}[column sep = 25, row sep = 18]
			\A \arrow[rr, bend right=60, "h"', ""{name=C}] \arrow[rr, "g" {near end, description}, ""{name=B}] \arrow[rr, bend left=60, "f", ""{name=A}]
			&& \B &&&& A_{0} \arrow[rr, "\left(\alpha{,}\beta\right)"] && B_{2} \arrow[rr, "m"] && B_{1}
			\arrow[Rightarrow, from=A, to=B, "\overline{\alpha}", shorten <= 4pt, xshift = -1.75ex]
			\arrow[Rightarrow, from=B, to=C, "\overline{\beta}", shorten <= 4pt, xshift = -1.75ex]
\end{tikzcd}$$ 

The left whiskering and right whiskering pictured below are defined as the composites in $\E$ given by $\beta f_{0}$ and $g_{1}\alpha$ respectively, and the horizontal composition of $2$-cells is defined via whiskering and vertical composition in the usual way as described in Proposition II 3.1 of \cite{mac2013categories}. 
 	
 	\begin{center}
 		\begin{tikzcd}[column sep = 12, row sep = 12]
 			\mathbb{A} \arrow[rr, "f"]
 			&& \mathbb{B} \arrow[rr, "g", bend left =50, ""{name=A}] \arrow[rr, "g'"', bend right=50, ""{name=B}]
 			&& \mathbb{C} &&{}
 			\arrow[Rightarrow, from=A, to=B, "\overline{\beta}", shorten <= 4pt]
 		\end{tikzcd} \begin{tikzcd}[column sep = 12, row sep = 12]
 			\mathbb{A} \arrow[rr, bend left=50, "f", ""{name=A}] \arrow[rr, bend right=50, "f'"', ""{name=B}]
 			&& \mathbb{B} \arrow[rr,"g"]
 			&& \mathbb{C}
 			\arrow[Rightarrow, from=A, to=B, "\overline{\alpha}", shorten <= 4pt]
 		\end{tikzcd} 
 	\end{center} 
  
   If $\mathcal{E}$ is small (resp. locally small), then $\mathbf{Cat}\left(\mathcal{E}\right)$ is small (resp. has small hom-categories).
 \end{prop}

 Further background on properties of the $2$-category $\CatE$ will be reviewed in Remark~\ref{Bourke black box}.

  Fully-faithfulness for internal functors is recalled in Definition~\ref{def:ff}, to follow. Unlike in the enriched setting, this is equivalent to the representably defined notion of fully-faithfulness for morphisms in $\CatE$.

 \begin{definition}\label{def:ff}
     Let $\E$ be a category with products. An internal functor $f: \A \to \B$ is called
     
     \begin{itemize}
         \item \emph{faithful} if the morphism into the pullback induced by the following commutative square is a monomorphism.
         \item \emph{fully faithful} if the induced morphism into the pullback is an isomorphism.
     \end{itemize} 

     \begin{equation*}
         \begin{tikzcd}
             A_1 \arrow[r, "f_1"] \arrow[d, "(d_0{,}d_1)", swap] & B_1 \arrow[d, "(d_0{,} d_1)"] \\
             A_0 \times A_0 \arrow[r, "f_0 \times f_0", swap] & B_0 \times B_0
         \end{tikzcd}
     \end{equation*}
 \end{definition}

\begin{remark}\label{foreshadowing remark full subobjects}
    An internal functor $(f_{0}, f_{1})$ is a monomorphism in $\CatE$ if and only if it is faithful and $f_{0}$ is a monomorphism. In Section~\ref{Section Subobject Classifiers} we will relate subobject classifiers in $\E$ to classifiers for morphisms in $\CatE$ which are both fully faithful and monomorphisms; these notions being definable representably in $\CatE$. In Subsection~\ref{AC in 2-categorical terms}, we will exhibit such morphisms in $\CatE$ as the right class $\R'$ of an orthogonal factorisation system, giving an internal version of the analysis in Section 5.2 of \cite{bourke2014two}. The left class $\Lagr'$ of this factorisation system will consist of internal functors $f: \A \to \B$ for which $f_{0}: A_{0} \to B_{0}$ are epimorphisms in $\E$. This will allow us to detect them via the $2$-category structure of $\CatE$, despite the fact that representables $\E(X, -): \E \to \s$ typically fail to preserve or jointly reflect epimorphisms. The class $\Lagr'$ features in our categorification of the axiom of choice, in Definition~\ref{definition categorified axiom of choice}.
	\end{remark}

\subsection{Adjunctions between $\E$ and $\mathbf{\CatE}_{1}$}\label{subsection adjunctions between E and Cat E}
  We review some adjunctions between $\CatE_1$ to $\E$. These adjunctions will be invaluable in our proofs that various universal properties in one of these categories imply analogous properties in the other. 

 \begin{remark}\label{Disc left adjoint counit}
	The functor $(-)_0: \CatE_1 \to \E$ has a left adjoint $\disc: \mathcal{E} \rightarrow \mathbf{Cat}\left(\mathcal{E}\right)_{1}$. This sends $X \in \mathcal{E}$ to the internal category $\disc(X): \Delta_{\leq 3}^\text{op} \to \E$ which is constant at $X$. The components of the unit of this adjunction on $X \in \mathcal{E}$ are all given by identities, and as such $\mathbf{disc}: \E \to \CatE_{1}$ is fully faithful. Indeed, it is the inclusion of the category of discrete objects in $\CatE$ in the sense of Remark~\ref{2-categorical conventions}. Meanwhile the components of the counit on an internal category $\mathbb{A}$ are given by the internal functor whose component on objects is $1_{A_{0}}$ and component on arrows is $i: A_{0} \to A_{1}$. It is easy to see that the naturality square for the counit on an internal functor $f$ is a pullback precisely if $f$ reflects identities, in the sense that the square $f_{1}i = if_{0}$ is a pullback. It is also easy to see that $\mathbf{disc}$ preserves finite limits, even when it does not have the left adjoint that will be described in Remark~\ref{Connected Components}.	
	\end{remark}

 \begin{remark}\label{CoDisc right adjoint}
     When $\E$ has products, $(-)_0: \CatE_{1} \to \E$ also has a right adjoint, which we call $\codisc: \E \to \CatE_1$. This sends $X$ to the internal category defined by $\{n \mapsto X^n\}$, with $n$-simplices given by the $n$-fold product for $n \in \Delta_{\leq 3}^\text{op}$. When $\E = \s$, this is the groupoid with set of objects is $X$ and a unique morphism between any two objects. The counit of $(-)_{0} \dashv \mathbf{indisc}$ is the identity, and as such $\mathbf{indisc}: \E \to \CatE_{1}$ is fully faithful. Meanwhile the unit has its component on an internal category $\mathbb{A}$ given by the internal functor $\eta_\mathbb{A}: \A \to \mathbf{indisc}(\A)$ which is given by the identity on objects, and the morphism $(d_{0}, d_{1}): A_{1} \to A_{0} \times A_{0}$ between objects of arrows. Observe that an internal functor $f$ is fully faithful if and only if the naturality square of $\eta^\mathcal{E}$ on $f$ is a pullback. 
     
      We call internal categories of the form $\codisc(X)$ for some $X \in \E$ \emph{indiscrete}. Note that there are other names for this in the literature: \emph{chaotic, codiscrete, coarse} and \emph{Brandt}.
     \end{remark}

     \begin{remark}\label{BO FF ofs on Cat E}
         The counit of $\mathbf{disc} \dashv (-)_{0}$ and the unit of $(-)_{0}\dashv \mathbf{indisc}$ both have components which are internal functors given by isomorphisms (indeed, identities) between objects of objects. Isomorphism on objects internal functors $f: \mathbb{A} \to \mathbb{B}$ play a special role in the $2$-category $\CatE$. They are strongly left orthogonal to fully faithful internal functors, in the sense of Definition 2.3.3 of \cite{bourke2010codescent}. Indeed, they form the left class of an orthogonal factorisation system in $\CatE_{1}$, for which the right class are the fully faithfuls. This factorisation is constructed via certain $2$-categorical limits and colimits in $\CatE$, which we will describe in more detail in Remark~\ref{Bourke black box}.
     \end{remark}

	\begin{remark}\label{Connected Components}
		Assume $\mathcal{E}$ has coequalisers of reflexive pairs. Then $\mathbf{disc}$ has a left adjoint $\Pi_{0}: \mathbf{Cat}\left(\mathcal{E}\right) \rightarrow \mathcal{E}$ which sends every internal category $\mathbb{A}$ to the codomain of the coequaliser $q_{A}$ of its source and target, and every internal functor $\left(f_{0}, f_{1}\right): \mathbb{A} \rightarrow \mathbb{B}$ to the morphism shown below, which is induced by the universal property of $\Pi_{0}(\mathbb{A})$, given the serial commutativity of the square on the left.
		
		\begin{center}
			\begin{tikzcd}[column sep = 15, row sep = 10]
				{A}_{1} \arrow[rr, shift left = 5pt, "{d}_{0}"] \arrow[rr, shift right = 5 pt, "{d}_{1}"'] \arrow[dd, "{f}_{1}"']
				&& {A}_{0} \arrow[rr, "{q}_\mathbb{A}"] \arrow[dd, "{f}_{0}"']
				&& \Pi_{0}(\A) \arrow[dd, dashed, "\Pi_{0}{\left(f\right)}"]
				\\
				\\
				{B}_{1} \arrow[rr, shift left = 5pt, "{d}_{0}"] \arrow[rr, shift right = 5pt, "{d}_{1}"']
				&& {B}_{0} \arrow[rr, "{q}_\mathbb{B}"']
				&& \Pi_{0}(\B)
			\end{tikzcd}
		\end{center}
		
		 Since $\disc: \E \to \CatE_{1}$ is fully faithful, the component of the counit on an object $X \in \mathcal{E}$ can again be chosen to be the identity. Meanwhile, the component of the unit $q: 1_{\mathbf{Cat}\left(\mathcal{E}\right)_{1}} \Rightarrow \mathbf{disc}\circ \Pi_{0}$ on an internal category $\mathbb{A}$ is given on objects by the coequaliser $q_\mathbb{A}$ above, and on arrows by the subsequent composite from $A_{1}$ to $\Pi_{0}(\A)$. The triangle identities can be shown using the universal properties of the coequalisers.
 
   For the proof of Theorem~\ref{Theorem generator in E vs generator in CatE} we will need the more nuanced observation that there is also a natural bijection $\E(\Pi_{0}(\A), B) \cong \CatE_{1}(\A, \mathbf{disc}(B))$, defined whenever the coequaliser of the source and target morphisms for $\A$ exists in $\E$. It is straightforward to see that this also holds, via a similar argument to the one sketched above.
	\end{remark}

\section{Finite limits and Bourke's characterisation of $\CatE$}\label{Section Finite Limits}

If $\E$ has pullbacks then on top of pullbacks, $\CatE$ also has powers by the category $\mathbf{2}$, containing the free-living arrow. These are given by an internal version of arrow categories, and will be described briefly in Remark~\ref{Bourke black box}. A more detailed explicit internal description is given in \cite{bourke2010codescent,miranda2022internal}. Moreover, $2$-categories of the form $\CatE$ have been characterised by Bourke, as we recall in Proposition~\ref{Bourke characterisation of Cat E} to follow. For our purposes, it suffices to know that $2$-categories of the form $\CatE$ may be characterised in elementary and purely $2$-categorical terms.

\begin{prop}\label{Bourke characterisation of Cat E}
    (Theorem 4.18 of \cite{bourke2010codescent}) If $\E$ is a category with pullbacks then the $2$-category $\mathcal{K}:= \CatE$ satisfies the conditions listed below. Conversely, if $\mathcal{K}$ satisfies the conditions listed below, then there is a $2$-equivalence $K \simeq \mathbf{Cat}\left(\E\right)$ where $\E:= \mathbf{Disc}\left(\mathcal{K}\right)$.

    \begin{enumerate}
        \item $\mathcal{K}$ has pullbacks and powers by $\mathbf{2}$.
        \item $\mathcal{K}$ has codescent objects of categories internal to $\mathcal{K}$ whose source and target maps form a two-sided discrete fibration.
        \item Codescent morphisms are effective in $\mathcal{K}$.
        \item Discrete objects in $\mathcal{K}$ are projective, in the sense of Definition 4.13 of \cite{bourke2010codescent}.
        \item For every object $A \in \mathcal{K}$, there is a projective object $P \in \mathcal{K}$ and a codescent morphism $c: P \rightarrow A$.
    \end{enumerate}
\end{prop}

\begin{remark}\label{Bourke black box}
    In this paper we mostly work with $2$-categories $\mathcal{K}$ which satisfy the conditions listed in Proposition~\ref{Bourke characterisation of Cat E}. When doing so, Bourke's result allows us to use the techniques of internal category theory in our proofs, even when dealing with properties stated in purely $2$-categorical terms.
    
     Although readers should be able to follow our proofs by treating Proposition~\ref{Bourke characterisation of Cat E} as a `black box', we give some brief comments on its content. The powers $\mathbb{A}^\mathbf{2}$ can be given an explicit description internally to $\E$. They have objects of objects given by $A_{1}$; the object of arrows of $\mathbb{A}$, while their objects of arrows are given by the pullback depicted below which may be thought of as the `object of internal squares in $\mathbb{A}$'.

    \begin{equation*}
        \begin{tikzcd}
            A_{Sq} \arrow[r] \arrow[d] \arrow[dr, phantom, "\lrcorner", very near start] & A_2 \arrow[d, "m"] \\
            A_2 \arrow[r, "m", swap] & A_1
        \end{tikzcd}
    \end{equation*}

 Codescent objects in a $2$-category $\mathcal{K}$ are $2$-categorical colimits of truncated simplicial objects, defined by the weight $\Delta_{\leq 2} \to \mathbf{Cat}$, where $\Delta_{\leq 2}$ is considered as a $2$-category with only identity $2$-cells. Categories internal to $\mathcal{K}$ whose source and target maps form a two-sided discrete fibration are called \emph{cateads} in $\mathcal{K}$. Such data can be thought of as two-dimensional versions of preorders, with the condition that $(d_{1}, d_{0}): A_{1} \to A_{0}\times A_{0}$ should be jointly monic being replaced by the condition that it should be a two-sided discrete opfibration.

 An object $A$ of a $2$-category $\mathcal{K}$ is said to be \emph{projective} if the representable $\mathcal{K}(\A, -): \mathcal{K} \to \Cat$ preserves codescent morphisms. This extends the one-dimensional notion, where instead the representable preserves regular epimorphisms. 

 Codescent morphisms for cateads in $\mathcal{K}= \CatE$ are precisely those internal functors $f: \mathbb{A} \to \mathbb{B}$ for which $f: A_{0} \to B_{0}$ are isomorphisms. One may think of a catead $\mathbb{C}$ in $\CatE$ as a two-dimensional version of an equivalence relation. From this perspective, its codescent object is a two-dimensional quotient, which is equivalently given by the `$0$-th row' of the underlying double category in $\E$. If the internal category of objects of the double category $\mathbb{C}$ is called its underlying vertical category internal to $\E$, then the codescent object of $\mathbb{C}$ is its underlying horizontal category internal to $\E$.

 As mentioned in Remark~\ref{BO FF ofs on Cat E}, (iso on objects, fully faithful) forms an orthogonal factorisation system on $\CatE_{1}$. We briefly review its construction. Given an internal functor $f: \mathbb{A} \to \mathbb{B}$, first form the following double category, or category internal to $\CatE_{1}$.

         $$\begin{tikzcd}
             f\downarrow f\downarrow f \arrow[rr, "m"description]\arrow[rr, shift left = 2,"\pi_{1}"]\arrow[rr, shift right = 2,"\pi_{0}"']&& f\downarrow f \arrow[rr, shift left = 2, "d_{1}"]\arrow[rr, shift right = 2, "d_{0}"'] && \mathbb{A}\arrow[ll, "i"description]
         \end{tikzcd}$$

          Where $f\downarrow f$ and $f \downarrow f \downarrow f$ are respectively given by the comma and pullback in $\CatE$ depicted below left and below right.

         $$\begin{tikzcd}
             f\downarrow f \arrow[d, "d_{0}"']\arrow[r, "d_{1}"] & \mathbb{A}\arrow[ld, Rightarrow, shorten = 10, "\phi"]\arrow[d, "f"]
             \\
             \mathbb{A}\arrow[r, "f"'] & \mathbb{B}&{}
         \end{tikzcd}\begin{tikzcd}
             f \downarrow f \downarrow f \arrow[dr, phantom, "\lrcorner", very near start]\arrow[r, "\pi_{0}"]\arrow[d, "\pi_{1}"'] & f\downarrow f\arrow[d, "d_{1}"]
             \\
             f\downarrow f\arrow[r, "d_{0}"'] & \mathbb{A}
         \end{tikzcd}$$

          Bourke shows that the double category just described is a catead, and that the factorisation $f =hk$ where $k$ is given by an isomorphism between objects of objects and $h$ is fully faithful, is given by taking $k: \mathbb{A} \to \mathbb{C}$ to be coprojection to the codescent object for this catead and $h: \mathbb{C} \to \mathbb{B}$ to be the internal functor induced by the universal property of $\mathbb{C}$. The adjective `effective' in part (3) of Proposition~\ref{Bourke characterisation of Cat E} then amounts to the fact that $h$ is an isomorphism in $\CatE$ if and only if $f_{0}: A_{0} \to B_{0}$ is an isomorphism in $\E$. Finally, projective covers are given by $\varepsilon_{\A}: \mathbf{disc}(\A)_{0} \to \A$; the components of the counit of the adjunction $\mathbf{disc} \dashv (-)_{0}$ described in Remark~\ref{Disc left adjoint counit}.
\end{remark}

\begin{prop}\label{prop finite weighted limits in Cat E}
    The $2$-category $\mathbf{Cat}\left(\E\right)$ has all finite $2$-limits if and only if the category $\E$ has all finite limits. 
\end{prop}

\begin{proof}
    By Proposition~\ref{Bourke characterisation of Cat E}, and the fact that $\mathbf{2}$ is a strong generator in $\Cat$, it suffices to show that $\E$ has terminal objects if and only if $\CatE$ does. But this follows from the adjunctions $\disc \dashv (-)_{0}\dashv \codisc$.
\end{proof}

\section{Cartesian closedness}\label{Section Cartesian Closedness}

Recall that exponentials $[X, Y]$ in $\s$ consist of sets whose elements are functions from $X$ to $Y$, while exponentials $[\mathcal{C}, \mathcal{D}]$ in $\Cat$ consist of categories whose objects are functors from $\mathcal{C}$ to $\mathcal{D}$, and whose morphisms are natural transformations between these functors. In this Section we consider an $\E$-internal version of these functor categories, which can also be constructed in terms of exponentials and finite limits in $\E$.

\begin{thm}\label{Theorem cartesian closedness iff}
    Let $\E$ be a category with finite limits. The category $\E$ is cartesian closed if and only if the $2$-category $\CatE$ is cartesian closed. In this case, $\mathbf{disc}: \E \to \CatE_{1}$ preserves internal homs.
\end{thm}

\begin{proof}
    Cartesian closedness of the category $\CatE_{1}$ has been shown in \cite{bastiani1972categories}, under the assumption that $\E$ has finite limits and exponentials, by viewing $\CatE_{1}$ as the category of models of a finite limit sketch. Indeed, it is shown in Theorem 2.1.1 of \cite{miranda2022internal} that the nerve $N: \CatE_{1} \to [\Delta^\text{op}, \E]$ is an inclusion of an exponential ideal. The two-dimensional aspect of the universal property of cartesian closedness for the $2$-category $\CatE$ follows from the universal property of powers by $\mathbf{2}$, which we denote as $\mathbf{2}\pitchfork (-)$. In particular, it is exhibited by the following natural bijections.
    
    $$\CatE_{1}(\mathbb{A} \times \mathbb{B}, \mathbf{2}\pitchfork{\mathbb{C}}) \cong \CatE_{1}(\mathbb{A} , \left({\mathbf{2}} \pitchfork{\mathbb{C}}\right)^{\mathbb{B}}) \cong \CatE_{1}(\mathbb{A} , \mathbf{2}\pitchfork\left({\mathbb{C}}^{\mathbb{B}}\right))$$

    Conversely, let $\E$ be a category with finite limits and suppose $\CatE$ is cartesian closed. We show that $\E$ is cartesian closed with exponentials given as displayed below for $Y, Z \in \E$. 
    
    $$Z^Y : = (\disc(Z)^{\disc(Y)})_0$$ 
    
    The following calculations show that the proposed exponential satisfies the isomorphism depicted below, naturally in all $X,Y,Z \in \E$.

    $$ \Hom(X \times Y, Z) \cong \Hom(X, Z^Y)$$

    \begin{align*}
        \E(X \times Y, Z) & = \E(X \times Y, \disc(Z)_0) &\text{ {(}unit of }\mathbf{disc}\dashv (-)_{0}\text{ is the identity{)}}\\
        & \cong \CatE_{1}(\disc(X \times Y), \disc(Z)) &{(}\mathbf{disc} \text{ is fully-faithful}{)} \\
        & \cong \CatE_{1}(\disc(X) \times \disc(Y), \disc(Z)) & {(}\mathbf{disc}\text{ preserves products}{)} \\
        & \cong \CatE_{1}\left(\disc(X) , \disc(Z)^{\disc(Y)}\right) &{(}\CatE\text{ is cartesian closed} {)}\\
        & \cong \E\left(X, \left(\disc(Z)^{\disc(Y)}\right)_0\right) & (\mathbf{disc} \dashv(-)_{0})\\
        & =: \E\left(X, Z^Y\right).&{}
    \end{align*}

    Cartesian closedness of $\mathbf{disc}: \E \to \CatE$ is an easy inspection given the construction of internal homs in $\CatE$, and also follows from Day's reflection theorem \cite{day1972reflection}.
\end{proof}

\section{Well-pointedness}\label{Section Generators}

Recall that in $\s$, we can test whether two functions $f,g: X \to Y$ are equal by checking if $f(x)=g(x)$ for every $x \in X$. Similarly, in $\Cat$, to test if two functors $F, G : \mathcal{C} \to \mathcal{D}$ are equal it suffices to check that $Ff = Gf$ for every $f \in \mathcal{C}_1$. This amounts to $\mathbf{1}$ being a generator for $\s$ and $\mathbf{2}$ being a generator for $\Cat$. The aim of this section is to show that the analogous statements for $\E$ and $\CatE$ are logically equivalent under the assumption that $\E$ is lextensive and cartesian closed. As we saw in Theorem~\ref{Theorem cartesian closedness iff}, $\E$ is cartesian closed if and only if $\CatE$ is. We first show a similar logical equivalence between extensivity of $\E$ and of $\CatE$. It will follow that $\E$ is lextensive if and only if $\CatE$ is.

\begin{definition}\label{Definition extensivity}
\begin{enumerate}
    \item A category with pullbacks $\E$ is said to be \emph{extensive} \cite{CARBONI1993145} if it has finite coproducts and for all $A, B \in \mathcal{E}$, the functor $\mathcal{E}/A \times \mathcal{E}/B \rightarrow \mathcal{E}/\left(A + B\right)$, which takes the coproduct, is an equivalence of categories. Call an extensive category \emph{lextensive} if it moreover has a terminal object.
    \item Call a $2$-category with pullbacks $\mathcal{K}$ extensive if it has finite coproducts and the similarly defined $2$-functor is a $2$-equivalence. Call an extensive $2$-category $\mathcal{K}$ lextensive if it moreover has a terminal object and powers by $\mathbf{2}$.
\end{enumerate}
    
\end{definition}

\begin{lem}\label{extensive coproduct Cat(E)}
    Let $\E$ be a category with pullbacks and products. The category $\mathcal{E}$ is extensive if and only if the $2$-category $\mathbf{Cat}\left(\mathcal{E}\right)$ is extensive, in which case the coproducts in $\mathbf{Cat}\left(\mathcal{E}\right)$ are computed in $[\Delta_{\leq 3}^\text{op}, \mathcal{E}]$.
\end{lem}

\begin{proof}
	It is clear from the adjunctions $\mathbf{disc}\dashv (-)_{0}\dashv \mathbf{indisc}$ that $\E$ has an initial object if and only if $\CatE_{1}$ does, and in this case so does the $2$-category $\CatE$. The functor category $[\Delta_{\leq 3}^\text{op}, \mathcal{E}]$ has whatever colimits $\mathcal{E}$ has, computed pointwise. Suppose $\E$ has extensive coproducts. Let $\mathbb{A}$ and $\mathbb{B}$ be categories internal to $\mathcal{E}$. Then the diagrams which need to be pullbacks for $\mathbb{A} + \mathbb{B}$ to be well-defined as an internal category are precisely the coproducts in $\mathcal{E}$ of the corresponding pullbacks which exhibit $\mathbb{A}$ and $\mathbb{B}$ as internal categories. But by extensivity of $\mathcal{E}$, these will be pullbacks as well. Thus the category $\mathbf{Cat}\left(\mathcal{E}\right)_{1}$ has coproducts as computed in $[\Delta_{\leq 3}^\text{op}, \mathcal{E}]$. But the two-dimensional aspect of the universal property for coproducts follows from the one-dimensional aspect, since $\CatE$ has powers by $\mathbf{2}$.

  Conversely, suppose that $\CatE$ has extensive coproducts. For $X, Y \in \E$, we claim that their coproduct is given by $(\mathbf{disc}(X) + \mathbf{disc}(Y))_{0}$. By Remark~\ref{CoDisc right adjoint}, since $\E$ has products the functor $(-)_{0}: \CatE_{1} \to \E$ is a left adjoint and hence preserves coproducts. But $(\mathbf{disc}(X))_{0} = X$ and $(\mathbf{disc}(Y))_{0} = Y$. This completes the proof.
\end{proof}

 For the remainder of this section we assume that $\E$, and hence $\CatE$, is lextensive.

Next, we recall the construction of the free-living arrow $\mathbf{2}_{\mathcal{E}}$ as a category internal to $\E$.  Copowers by $\mathbf{2}$ in $\CatE$ can be constructed in terms of this internal category, as we will show in Theorem~\ref{thm:copower}.

\begin{remark}
\label{remark construction F : Setf -> E}
    Recall that any finite limit preserving functor between finite limit categories $G: \mathcal{S} \rightarrow \mathcal{E}$ gives rise to a $2$-functor $\mathbf{Cat}\left(G\right): \mathbf{Cat}\left(\mathcal{S}\right) \rightarrow \mathbf{Cat}\left(\mathcal{E}\right)$, which acts componentwisely on all data \cite{miranda2022internal}. Recall also that the category of finite sets $\mathbf{FinSet}$ is the free completion under finite coproducts of the terminal category. Furthermore, for lextensive $\mathcal{E}$, the unique coproduct preserving functor $F_\mathcal{E}: \mathbf{FinSet} \rightarrow \mathcal{E}$ which preserves the terminal object also preserves all other finite limits. 
\end{remark}

\begin{definition}\label{def free living arrow}
    Take $\mathcal{S} = \mathbf{FinSet}$ as in Remark~\ref{remark construction F : Setf -> E} and apply the $2$-functor $\mathbf{Cat}\left(F_\mathcal{E}\right): \mathbf{Cat}\left(\mathbf{FinSet}\right) \rightarrow \mathbf{Cat}\left(\mathcal{E}\right)$ to the free living arrow $\mathbf{2} \in \mathbf{Cat}\left(\mathbf{FinSet}\right)$. Denote the resulting category internal to $\E$ as $\mathbf{2}_\mathcal{E}$. 
\end{definition}

 The internal category $\mathbf{2}_\mathcal{E}$ of Definition~\ref{def free living arrow} can be described explicitly as a truncated simplicial object, with $n$-simplices given by the $(n+2)$-fold coproduct of the terminal object $\mathbf{1} \in \E$; see Example 2.3.2 of \cite{miranda2022internal} for details. Recall that the copower by $\mathbf{2}$ of an object $A \in \mathcal{K}$, if it exists, is an object $ \mathbf{2} \odot A$ equipped with isomorphisms of categories $\mathcal{K}\left(\mathbf{2} \odot A, B\right) \cong \mathbf{Cat}\left(\mathbf{2}, \mathcal{K}\left(A, B\right)\right)$ which vary $2$-naturally in $B$. The next theorem then shows that the $2$-functor $\mathbf{Cat}\left(F_\mathcal{E}\right): \mathbf{Cat}\left(\mathbf{FinSet}\right) \rightarrow \mathbf{Cat}\left(\mathcal{E}\right)$ preserves copowers by $\mathbf{2}$.

\begin{thm}\label{thm:copower}
	Let $\mathcal{E}$ be lextensive and cartesian closed, and let $\mathbf{2}_\E$ be constructed as in Definition~\ref{def free living arrow}.
 
 \begin{enumerate}
    \item The internal hom $[\mathbf{2}_\mathcal{E}, \mathbb{B}]$ has the universal property of the power of $\mathbb{B}$ by $\mathbf{2}$. 
     \item For $\mathbb{A} \in \mathbf{Cat}\left(\mathcal{E}\right)$, the internal category $\mathbf{2}_\mathcal{E} \times \mathbb{A}$ has the universal property of the copower of $\mathbb{A}$ by $\mathbf{2}$ in $\mathbf{Cat}\left(\mathcal{E}\right)$.
 \end{enumerate}
\end{thm}

\begin{proof}
	 Consider the unique non-identity naturnal transformation $\rho$ from the category $\mathbf{1}$ to the category $\mathbf{2}$. The internal functor $[\mathbf{2}_\mathcal{E}, \mathbb{B}] \rightarrow \mathbb{B}^\mathbf{2}$ is induced by the universal property of the power by $\mathbf{2}$ given the image of $\rho$ under the $2$-functor displayed below. 
  
  $$\begin{tikzcd}
	     \mathbf{Cat}(\mathbf{FinSet})^\text{op} \arrow[rr, "\Cat(F_\E)^\text{op}"] && \CatE^\text{op} \arrow[rr, "{[}-{,}\B{]}"] &&\CatE
	 \end{tikzcd}$$
  
  We describe the transpose $\mathbf{2}_\mathcal{E} \times \mathbb{B}^\mathbf{2} \rightarrow \mathbb{B}$ of the required inverse internal functor $\mathbb{B}^\mathbf{2} \rightarrow [\mathbf{2}_\mathcal{E}, \mathbb{B}]$. Recall first that $\mathbf{2}$ is the category that has, as objects, the set $\{*\} + \{*\}$ and, as arrows, the set $\{*\} + \{*\} + \{*\}$. By lextensivity of $\E$ and as $\Cat(F_{\E})$ preserves coproducts, $\mathbf{2}_\E \times \B$ has, as objects, $B_1 + B_1$ and as arrows $B_{\text{sq}}+ B_{\text{sq}} + B_{\text{sq}}$. Now, between objects of objects, the functor $\mathbf{2}_\mathcal{E} \times \mathbb{B}^\mathbf{2} \rightarrow \mathbb{B}$ is given by  $\left(d_{0}, d_{1}\right): B_{1} + B_{1} \rightarrow B_{0}$ induced by the universal property of the coproduct, using the source and target maps. Between objects of arrows it is given by the morphism $B_\text{sq} +B_\text{sq} +B_\text{sq} \rightarrow B_{1}$ induced by the universal property of the coproduct by the source and target maps of $\mathbb{B}^\mathbf{2}$, as well as by the diagonal of the pullback square defining $B_\text{sq}$. To prove internal functoriality, one needs to check commutativity conditions for maps out of coproducts. These can in turn be verified by checking cases for each summand appearing in the coproduct. However, each of these individual cases just involves pullbacks and hence follows from the analogous property when $\mathcal{E}=\mathbf{Set}$, using the Yoneda Lemma. The proof that these internal functors are mutually inverse is similar. This proves part (1). Part (2) then follows by the following chain of isomorphisms, where the penultimate step uses part (1).
\begin{align*}
    \mathbf{Cat}\left(\mathbf{2}, \mathbf{Cat}\left(\mathcal{E}\right)\left(\mathbb{A}, \mathbb{B}\right)\right)&\cong \mathbf{Cat}\left(\mathcal{E}\right)\left(\mathbb{A}, \mathbb{B}\right)^\mathbf{2}
    \\
    &\cong \mathbf{Cat}\left(\mathcal{E}\right)\left(\mathbb{A}, \mathbb{B}^\mathbf{2}\right)
    \\
    &\cong \mathbf{Cat}\left(\mathcal{E}\right)\left(\mathbb{A}, [\mathbf{2}_\mathcal{E}, \mathbb{B}]\right)
    \\
    &\cong \mathbf{Cat}\left(\mathcal{E}\right)\left(\mathbf{2}_\mathcal{E}\times\mathbb{A}, \mathbb{B}\right)
\end{align*}

\end{proof}

 In particular we have that $\mathbf{2}_{\E}$ has the universal property of the copower by $\mathbf{2}$ of the terminal object in $\CatE$.

\begin{remark}
    The assumptions of Theorem~\ref{thm:copower} part (2) can be relaxed. In particular, lextensivity of $\E$ suffices for $\CatE$ to have copowers by $\mathbf{2}$. One can directly check that $\mathbf{2}_{\E} \times \mathbb{A}$ has the appropriate universal property. However, doing so requires significant tedious calculations. Some of these calculations can be found in the Appendix of \cite{miranda2022internal}. We do not need this extra level of generality however since for our purposes we may assume that the category $\E$, or equivalently, the $2$-category $\CatE$, is cartesian closed. 
\end{remark}

 \begin{remark}\label{discrete natural transformations}
	Generating families, in the sense we will recall in Definition~\ref{definition generator}, can be constructed in $\CatE$ using copowers by $\mathbf{2}$. To show this we will need to observe that internal natural transformations out of discrete categories correspond to morphisms into the object of arrows of their codomain internal category. We now explain why this is so.
 
  Let $X \in \mathcal{E}$ and $\mathbb{A} \in \mathbf{Cat}\left(\mathcal{E}\right)$, and recall the adjunction $\mathbf{disc} \dashv (-)_0$ from Remark~\ref{Disc left adjoint counit}. Then there are the following natural bijections:
	
	$$\mathcal{E}\left(X, A_{1}\right) = \mathcal{E}\left(X, \left(\mathbb{A}^\mathbf{2}\right)_0\right) \cong \mathbf{Cat}(\E)_{1}\big(\mathbf{disc}\left(X\right), \mathbb{A}^\mathbf{2}\big) \cong [\mathbf{2}, \CatE\left(\mathbf{disc}\left(X\right), \mathbb{A}\right)]_{0}$$
	
	 Thus morphisms from $X$ to the object of arrows of an internal category $\mathbb{A}$ are in natural bijection with internal natural transformations between internal functors from the discrete category on $X$ to $\mathbb{A}$.
\end{remark}

\begin{definition}\label{definition generator}

    \begin{enumerate}
        \item A family of objects $\mathcal{G}$ in a category $\mathcal{C}$ is said to be \emph{generating} if the family of hom-functors $\mathcal{C}\left(X, -\right): \mathcal{C} \to \s$ for $X \in \mathcal{G}$ are jointly faithful.
        \item A family of objects $\widehat{\mathcal{G}}$ in a $2$-category $\mathcal{K}$ is said to be \emph{generating} if the family of hom-functors $\mathcal{K}\left(X, -\right): \mathcal{K} \to \Cat$ for $X \in \mathcal{G}$ are jointly faithful on $1$-cells and $2$-cells.
    \end{enumerate}
\end{definition}

\begin{cor}\label{cor:generator}
	 Suppose that $\mathcal{E}$ has finite limits, extensive coproducts, and a generating family of objects $\mathcal{G}$. Form the family of internal categories $\widehat{\mathcal{G}} := \{\mathbf{2}_\mathcal{E} \times \mathbf{disc}\left(X\right)| X \in \mathcal{G}\}$. Then $\widehat{\mathcal{G}}$ is a generating family for $\mathbf{Cat}\left(\mathcal{E}\right)$.
\end{cor}

\begin{proof}
	Let $f, g: \mathbb{A} \rightarrow \mathbb{B}$ be internal functors and assume that $fh = gh$ for all internal functors $h: \mathbf{2}_\mathcal{E} \times \mathbf{disc}\left(X\right) \rightarrow \mathbb{A}$ where $X \in \mathcal{G}$. By Proposition~\ref{arr faithful} part (1), to show that $\widehat{\mathcal{G}}$ is a generating family, it suffices to show that $f_{1} = g_{1}$ under this assumption. Denote by $\alpha: X \to A_{1}$ the component assigner of the internal natural transformation which corresponds to $h$ via the universal property of the copower by $\mathbf{2}$. Then the whiskerings $f\overline{\alpha} = g\overline{\alpha}$ are also equal in $\CatE$. But by Remark~\ref{discrete natural transformations}, any morphism $X \rightarrow A_{1}$ is $\mathcal{E}$ corresponds to an internal natural transformation between internal functors from $\mathbf{disc}\left(X\right)$ to $\mathbb{A}$ This amounts to saying that $f_{1}\alpha = g_{1}\alpha$ for all $\alpha: X \rightarrow A_{1}$, and hence $f_{1} = g_{1}$ as $X \in \mathcal{G}$.
 
  This shows that the family of $2$-functors $\CatE(G, -): \CatE \to \Cat$ for $G \in \widehat{\mathcal{G}}$ are jointly faithful on $1$-cells. But joint faithfulness on $2$-cells follows from joint faithfulness on $1$-cells as $\CatE$ has powers by $\mathbf{2}$. A parallel pair of internal natural transformations as depicted below left corresponds to a parallel pair of internal functors as depicted below right.

$$\begin{tikzcd}
    \A \arrow[rr, bend left = 30, "f"name=A]\arrow[rr, bend right = 30, "g"'name=B] &&\B&{}
    \arrow[from=A, to=B, Rightarrow, shift left = 2, "\overline{\beta}", shorten = 1mm]\arrow[from=A, to=B, Rightarrow, shift right = 2, "\overline{\gamma}"', shorten = 1mm]
\end{tikzcd}\begin{tikzcd}
&&\B
\\
    \A \arrow[rru, "f"]\arrow[rrd, "g"']\arrow[rr, shift left = 1, "\tilde{\beta}"]\arrow[rr, shift right = 1,"\tilde{\gamma}"']&& \B^\mathbf{2}\arrow[u, "d_{1}"']\arrow[d, "d_{0}"]
  \\
  &&\B
\end{tikzcd}$$

 By the one-dimensional aspect of $\widehat{\mathcal{G}}$ being a generator, the equality of such a pair of internal functors can be detected via $\CatE(G, -): \CatE \to \Cat$. As such, the equality of the original parallel pair internal natural transformations can also be detected via these representables.
 
\end{proof}

\begin{example}\label{example generators in presheaf categories}
    Let $\mathcal{C}$ be a small category and $\E := [\mathcal{C}^\text{op}, \mathbf{Set}]$. Then $\mathcal{E}$ has a generating family given by the representables $\mathcal{G}:= \{\mathcal{C}(-, X): \mathcal{C}^\text{op} \to \s| X \in \mathcal{C}\}$. Now, the $2$-functor $\mathbf{Cat}(-): \mathbf{LEX} \to 2\text{-}\mathbf{CAT}$ of Proposition 3.1.5 in \cite{miranda2022internal} preserves powers by small categories, and as such there is an isomorphism of $2$-categories $\mathbf{Cat}([\mathcal{C}^\text{op}, \s]) \cong [\mathcal{C}^\text{op}, \mathbf{Cat}]$, where the second of these is the $\mathcal{V}$-enriched functor category with $\mathcal{V}=\Cat$ and $\mathcal{C}$ considered as a $2$-category with only identity $2$-cells. Since $[\mathcal{C}^\text{op}, \mathbf{Cat}]$ is an enriched fuctor category, it has copowers computed pointwisely in $\mathbf{Cat}$. Corollary~\ref{cor:generator} then says that the following is a generating family for the $2$-category $[\mathcal{C}^\text{op}, \mathbf{Cat}]$. This coincides with the generating family for $[\mathcal{C}^\text{op}, \Cat]$ in terms of representables and copowers by the strong generator $\{\mathbf{2}\} \subseteq \mathbf{Cat}$.

    $$\left\{\ Y \mapsto \coprod_{f\in \mathcal{C}(Y{,}X)}\mathbf{2} \qquad\middle| \qquad X \in \mathcal{C}\right\}$$
\end{example}

\begin{remark}\label{Remark grenerators converse}
    If certain colimits exist in $\E$, then a generating family $\mathcal{G}'\subseteq \CatE$ also gives rise to a generating family on $\E$. Specifically, we need $\E$ to have coequalisers for all reflexive pairs of source and target morphisms where $\mathbb{G} \in \mathcal{G}'$. 
    
    $$\begin{tikzcd}
        G_{1} \arrow[rr,shift left = 2, "d_{0}"]\arrow[rr, shift right =2, "d_{1}"']&& G_{0}\arrow[rr, "q_{\mathbb{G}}"] && \Pi_{0}(\mathbb{G})
    \end{tikzcd}$$ 
    
     In this case, the partial adjunction $\E(\Pi_{0}(\mathbb{G}), X) \cong \CatE_{1}(\mathbb{G}, \mathbf{disc}(X))$ exists for all $\mathbb{G} \in \mathcal{G}$. The generating family in $\E$ is then given by $\mathcal{G} := \{\Pi_{0}(\mathbb{G})|\mathbb{G} \in \mathcal{G}\}$. We give a detailed proof only of a special case in Theorem~\ref{Theorem generator in E vs generator in CatE}, since this will be enough for our main results and since generating families are in practice typically easier to construct in $\E$ than in $\CatE$. The proof of this special case requires no extra colimit assumptions on $\E$. We leave the straightforward generalisation to the setting described here to the interested reader.
\end{remark}

 Recall that a category $\E$ is called \emph{well-pointed} if it has a terminal object $\mathbf{1}$ and the family containing just $\mathbf{1}\in \E$ is a generator. We introduce the following categorified version of this definition.

\begin{definition}\label{def:2-well-pointed}
    A $2$-category $\mathcal{K}$ is called \emph{$2$-well-pointed} if the following conditions hold.
    
    \begin{enumerate}
        \item $\mathcal{K}$ has a terminal object $\underline{\mathbf{1}}$.
        \item The copower $\mathbf{2}\odot \underline{\mathbf{1}}$ exists in $\mathcal{K}$.
        \item The family containing just $\mathbf{2}\odot \underline{\mathbf{1}}$ is a generator for $\mathcal{K}$, in the sense of Definition~\ref{definition generator} part (2).
    \end{enumerate}
\end{definition}

 There is one final lemma that we will need before we are ready to prove the main result of this section.

\begin{lem}\label{equaliser is diagonal}
    Let $\mathcal{C}$ be a category with finite products. For $A \in \mathcal{C}$, consider the diagram displayed below in which the morphisms $\Delta_{A}: A \to A \times A$ denotes the diagonal $(1_{A}, 1_{A})$ and \begin{tikzcd}
        A 
        &A\times A\arrow[l, "\pi_{1}"']\arrow[r, "\pi_{2}"]
        &A
    \end{tikzcd} denote the product projections. This diagram is an equaliser.

    $$\begin{tikzcd}
        A \arrow[rr, "\Delta_{A}"] && A \times A \arrow[rr,shift left = 2, "\pi_{1}\times \Delta_{A}"]\arrow[rr,shift right = 2, "\Delta_{A}\times \pi_{2}"'] && A \times A\times A
    \end{tikzcd}$$
\end{lem}

\begin{proof}
    This is straightforward to check when $\mathcal{C} = \mathbf{Set}$: the functions being equalised send $(x, y)$ to $(x, x, y)$ and $(x, y, y)$ respectively. These outputs are indeed equal precisely when $x=y$. The claim then follows representably for a general $\mathcal{C}$ with finite limits.
\end{proof}

\begin{thm}\label{Theorem generator in E vs generator in CatE}
    Let $\E$ be a lextensive, cartesian closed category. Then $\E$ is well-pointed if and only if $\CatE$ is $2$-well-pointed in the sense of Definition~\ref{def:2-well-pointed}.
\end{thm}

\begin{proof}
    Recall that by Theorem~\ref{thm:copower}, the copower $\mathbf{2}\odot \underline{\mathbf{1}} \in \CatE$ may be taken as $\mathbf{2}_\E$. Corollary~\ref{cor:generator} therefore specialises to show that $\E$ being well-pointed implies that $\CatE$ is $2$-well-pointed by taking $\mathcal{G}:= \{\mathbf{1}\}$. For the converse, recall from Remark~\ref{Connected Components} that a natural bijection $\E(\Pi_{0}(\A), B) \cong \CatE_{1}(\A, \mathbf{disc}(B))$ exists if the source and target morphisms of $\A$ have a coequaliser in $\E$. Recall from the discussion after Definition~\ref{def free living arrow}, with further details found in Example 2.3.2 of \cite{miranda2022internal}, that the internal category $2_{\E}$ has object of $n$-simplices given by the $(n+2)$-fold coproduct of the terminal object $\mathbf{1}$. Now, Lemma~\ref{equaliser is diagonal} applies to $\mathcal{C}:= \E^\text{op}$ with $A = \mathbf{1}$, and shows that this coequaliser does exist in $\E$, so that $\Pi_0(\mathbf{2}_{\E}) \cong \mathbf{1}$. Therefore 
 $\E(\Pi_{0}(\mathbf{2}_\E), B) \cong \CatE_{1}(\mathbf{2}_{\E}, \mathbf{disc}(B))$ and so there is a bijection between diagrams of the following forms for $f, g \in \E(X, Y)$.

\begin{equation*}
\begin{tikzcd}
     \mathbf{1} \arrow[r] & X \arrow[r, "f", shift left= 2] \arrow[r, "g", swap, shift right] & Y
\end{tikzcd} 
\qquad
\begin{tikzcd}
     \mathbf{2}_{\E} \arrow[r] & \disc(X) \arrow[r, "\disc(f)", shift left= 2] \arrow[r, "\disc(g)", swap, shift right] & \disc(Y)
\end{tikzcd} 
\end{equation*}

 Hence if $\{\mathbf{2}_{\E}\}$ is a generator in $\CatE$ then $\mathbf{1}$ is a generator in $\E$. This completes the proof.
\end{proof}

 Observe that the assumptions of Theorem~\ref{Theorem generator in E vs generator in CatE} hold if $\E$ is an elementary topos. In Theorem~\ref{Theorem subobject classifier in E iff full subobject classifier in Cat E} we will characterise this stronger property for $\E$ in terms of $\CatE$. Observe also that copowers by $\mathbf{2}$ in $\CatE$ exist under assumptions which have already been shown to be equivalent for $\E$ and $\CatE$, namely lextensivity and cartesian closedness. Since $\CatE$ has copowers by $\mathbf{2}$, two-dimensional aspects of universal properties for $2$-limits can be inferred from the one-dimensional aspects of these universal properties. This is dual to the argument for faithfulness on $2$-cells of the family of $2$-functors $\CatE(G, -): \CatE \to \Cat$ for $G \in \widehat{\mathcal{G}}$, given in the proof of Corollary~\ref{cor:generator}. As such we will herein omit verification of two-dimensional aspects of universal properties for limits.

\section{Natural numbers objects}\label{Section NNO}

We show that $\E$ has a natural numbers object if and only if $\CatE$ has a natural numbers object, in the sense of Definition~\ref{definition NNO}, to follow. In particular, the work of this section shows that a natural numbers object in $\CatE$ is discrete on the natural numbers object of $\E$. Throughout this section we assume only that $\E$ has finite limits.

\begin{definition}\label{definition NNO}

    \begin{enumerate}
        \item Let $\mathcal{C}$ be a category with a terminal object $\mathbf{1}$. The data \begin{tikzcd}
            \mathbf{1} \arrow[r, "z"]& N \arrow[r, "s"] & N
        \end{tikzcd} is called a \emph{natural numbers object} in $\mathcal{C}$ if for any \begin{tikzcd}
            \mathbf{1} \arrow[r, "f"]& X \arrow[r, "g"] & X
        \end{tikzcd} there is a unique $u: N \to X$ making the diagram below commute.

        $$\begin{tikzcd}
            \mathbf{1}\arrow[rd, "f"', bend right] \arrow[r, "z"] & N \arrow[r, "s"]\arrow[d, dashed, "u"]&N\arrow[d, dashed, "u"]
            \\
            &X\arrow[r, "g"']&X
        \end{tikzcd}$$

        \item Let $\mathcal{K}$ be a $2$-category with a terminal object $\underline{\mathbf{1}}$. The data \begin{tikzcd}
            \underline{\mathbf{1}} \arrow[r, "z"]& \underline{N} \arrow[r, "s"] & \underline{N}
        \end{tikzcd} is called a \emph{natural numbers object} in $\mathcal{K}$ if it is a natural numbers object for the underlying $1$-category of $\mathcal{K}$ and, additionally, if given \begin{tikzcd}
            \underline{\mathbf{1}} \arrow[r, "f"]& X \arrow[r, "g"] & X
        \end{tikzcd} and \begin{tikzcd}
            \underline{\mathbf{1}} \arrow[r, "f'"]& X \arrow[r, "g"] & X
        \end{tikzcd} 
        which have corresponding maps $u,u': \underline{N} \to X$ respectively, whenever we have a $2$-cell $$\begin{tikzcd}
            \underline{\mathbf{1}} \arrow[r, "f" name=A, bend left] \arrow[r, "f'"' name=B, bend right] & X \arrow[from=A, to=B, Rightarrow, "\overline{\alpha}", shorten = 1mm]
        \end{tikzcd}$$ 

        then there is a unique $2$-cell as depicted below left, making the pasting diagram depicted below right commute.

        $$\begin{tikzcd}
            \underline{N} \arrow[r, bend left = 30, "u"name=A]\arrow[r, bend right = 30, "u'"'name=B] & X&{}\arrow[from=A, to=B, Rightarrow, "\overline{\phi}", shorten = 1mm]
        \end{tikzcd} \begin{tikzcd}[row sep = large]
            \underline{\mathbf{1}}\arrow[rrd, "f" name = X] \arrow[rrd, "f'"' name = Y, bend right= 50]\arrow[rr, "z"] && \underline{N} \arrow[rr, "s"]\arrow[d, bend left = 40, "u"name = A]\arrow[d,bend right = 40, "u'"'name=B]
            &&\underline{N}\arrow[d,bend left = 40, "u"name=C]\arrow[d,bend right = 40,"u'"'name=D]
            \\
            &&X\arrow[rr, "g"']
            &&X\arrow[from=A, to=B, Rightarrow, dashed, "\overline{\phi}", shorten = 1mm]\arrow[from=C, to=D, Rightarrow, dashed, "\overline{\phi}", shorten = 1mm] \arrow[from = X, to=Y, Rightarrow, "\overline{\alpha}", shorten = 1mm]
        \end{tikzcd}$$
        \end{enumerate}
\end{definition}

 We first give a proof of the following standard result. 

\begin{lem}
\label{lem: NNO preserved by LA}
    Let $\mathcal{C}, \mathcal{D}$ be categories with a terminal object and suppose $\mathcal{D}$ has a natural numbers object $(N, z: \mathbf{1} \to N, s: N \to N)$. If $L: \mathcal{D} \to \mathcal{C}$ is a left adjoint such that the unique morphism $j: L\mathbf{1} \to \mathbf{1}$ is invertible, then $(LN, L(z)\circ j^{-1}: \mathbf{1} \to LN, Ls: LN \to LN)$ is a natural numbers object for $\mathcal{C}$.
\end{lem}

\begin{proof}
    Let $R: \mathcal{C} \to \mathcal{D}$ be the right adjoint of $L$. As a right adjoint, the unique morphism $k: R\mathbf{1} \to \mathbf{1}$ is invertible. Let \begin{tikzcd} \mathbf{1} \arrow[r, "f"]& X \arrow[r, "g"] & X \end{tikzcd} be in $\mathcal{C}$. By the adjunction $L \dashv R$, there is a bijection between diagrams of the forms depicted below.

    \begin{equation*}
        \begin{tikzcd}
            \mathbf{1} \arrow[r, "z"] \arrow[dr, bend right, "R(z')\circ k^{-1}", swap] 
            & N \arrow[r, "s"] \arrow[d, dashed, "v"] 
            & N \arrow[d, dashed, "v"] \\
            & RX \arrow[r, "Rg", swap] & RX
        \end{tikzcd}
        \qquad
        \begin{tikzcd}[column sep = 35]
            \mathbf{1} \arrow[r, "L(z)\circ j^{-1}"] \arrow[dr, bend right, "f", swap] & LN \arrow[r, "Ls"] \arrow[d,dashed, "u"] & LN \arrow[d, dashed,"u"] \\
            & X \arrow[r, "g", swap] & X
        \end{tikzcd}
    \end{equation*}

     By the universal property of the natural numbers object $(N, z: \mathbf{1} \to N, s: N \to N)$ in $\mathcal{D},$ there is a unique such $v: N \to RX$ . Hence such a $u: LN \to X$ exists and is unique, as required. 
\end{proof}

 We obtain the following for one-dimensional natural number objects.

\begin{cor}
\label{cor NNO one dimensional iff}
    Let $\E$ be a category with terminal object and pullbacks. Then $\E$ has a natural numbers object if and only if $\CatE_1$ has a natural numbers object. 
\end{cor}

\begin{proof}
    Apply Lemma~\ref{lem: NNO preserved by LA} to $\disc(-) \dashv (-)_0$ for one implication, and to $(-)_0 \dashv \codisc(-)$ for the converse.
\end{proof}

 We extend this to a correspondence between a one-dimensional natural numbers object of $\E$ and a two-dimensional natural numbers object for $\CatE.$

\begin{thm}\label{thm NNO iff}
    Let $\E$ be a category with finite limits. Then $\E$ has a natural numbers object if and only if the $2$-category $\CatE$ has a natural numbers object. In this case, the functors $\mathbf{disc}: \E \to \CatE_{1}$, $(-)_{0}: \CatE_{1} \to \E$ and $\Pi_{0}: \CatE_{1} \to \E$ all preserve the natural numbers object.
\end{thm}

\begin{proof}
    By Corollary~\ref{cor NNO one dimensional iff}, $\E$ has a natural numbers object if and only if $\CatE_{1}$ does, and by Lemma~\ref{lem: NNO preserved by LA} the functors mentioned preserve the natural numbers object. It suffices to show that $(\disc(N), \disc(z):\underline{\mathbf{1}} \to \disc(N), \disc(s): \disc(N) \to \disc(N))$ satisfies the two-dimensional aspect of the universal property in Definition~\ref{definition NNO} part (2).
    
    Consider a diagram 
    
    $$\begin{tikzcd}
        \mathbf{1} \arrow[r, shift left = 2, "f"name = A, bend left]\arrow[r, shift right = 2, "f'"'name=B, bend right] & \mathbb{X} \arrow[r, "g"] &\mathbb{X}\arrow[from=A, to=B, Rightarrow, shorten = 3, "\overline{\alpha}"]
    \end{tikzcd}$$
    
    in $\CatE$ and let $u, u': \mathbf{disc}(N) \to \mathbb{X}$ be the morphisms induced by the universal property of the natural numbers object. Since $\underline{\mathbf{1}} = \disc(\mathbf{1}),$ the internal natural transformation $\overline{\alpha}: f \Rightarrow f'$ is uniquely determined by a map $\alpha: \mathbf{1} \to X_1$ in $\E$. This along with the morphism $g_{1}: X_{1} \to X_{1}$ uniquely determines a map $\phi: N \to X_{1}$ giving rise to an internal natural transformation satisfying the commutativity conditions in $\CatE$ depicted below left. Conversely, using the universal property of $\mathbb{X}^\mathbf{2}$ and the fact that $(\mathbb{X}^\mathbf{2})_{0} = X_1$, an internal natural transformation $\overline{\phi}: u\Rightarrow u'$ satisfying $g.\overline{\phi} = \overline{\phi}.s$ corresponds to a morphism in $\E$ satisfying the commutativity condition depicted below right, where $\alpha := \phi.z$. Hence, $\overline{\phi}$ is unique and $(\disc(N), \disc(z):\underline{\mathbf{1}} \to \disc(N), \disc(s): \disc(N) \to \disc(N))$ is a natural numbers objects for the $2$-category $\CatE$.

    $$ \begin{tikzcd}[row sep = large]
            \mathbf{1}\arrow[rrd, "f" name = X] \arrow[rrd, "f'"' name = Y, bend right= 50]\arrow[rr, "z"] && \disc(N) \arrow[rr, "s"]\arrow[d, bend left = 40, "u"name = A]\arrow[d,bend right = 40, "u'"'name=B]
            && \disc(N) \arrow[d,bend left = 40, "u"name=C]\arrow[d,bend right = 40,"u'"'name=D]
            \\
            &&\mathbb{X}\arrow[rr, "g"']
            &&\mathbb{X}\arrow[from=A, to=B, Rightarrow, dashed, "\overline{\phi}", shorten = 1mm]\arrow[from=C, to=D, Rightarrow, dashed, "\overline{\phi}", shorten = 1mm] \arrow[from = X, to=Y, Rightarrow, "\overline{\alpha}", shorten = 1mm]
        \end{tikzcd}
        \qquad
        \begin{tikzcd}
            \mathbf{1}\arrow[rd, "\alpha"', bend right] \arrow[r, "z"] & N \arrow[r, "s"]\arrow[d, dashed, "\phi"]&N\arrow[d, dashed, "\phi"]
            \\
            &X_1\arrow[r, "g_1"']&X_1
        \end{tikzcd}    
        $$
\end{proof}

\begin{remark}\label{remark coequalisers need NNOs}
    The category $\CatE_{1}$ may fail to have coequalisers even if $\E$ is an elementary topos. For example, take $\E:= \mathbf{FinSet}$, the category of finite sets. Then the parallel pair in $\CatE_{1}$ displayed below does not have a coequaliser. 
    
    $$\begin{tikzcd}
        \mathbf{1} \arrow[rr, shift left = 2, "d_{1}"]\arrow[rr, shift right = 2, "d_{0}"'] && \mathbf{2}
    \end{tikzcd}$$ 
    
     Indeed, the coequaliser of this parallel pair in $\Cat = \mathbf{Cat}(\s)$ is the monoid of natural numbers, considered as a one object category. As such, natural numbers objects seem to be necessary for the category $\CatE_{1}$ to have coequalisers, and hence for the $2$-category $\CatE$ to have finite $2$-colimits. Indeed, Lawvere observed in \cite{lawvere1966category} that coequalisers of functors between categories implies the `axiom of infinity'. This complexity remains for coinserters, even though they are PIE colimits \cite{power1991characterization}; the coinserter of the parallel pair in $\mathbf{Cat}$ is again the monoid of natural numbers. We will comment further on coequalisers in $\CatE$ in the conclusion, but leave detailed investigation to future research.
     
     $$\begin{tikzcd}
        \mathbf{1}\arrow[r,shift left =2, "1_\mathbf{1}"]\arrow[r, shift right = 2,"1_\mathbf{1}"']
        &\mathbf{1}
    \end{tikzcd}$$ 
\end{remark}

\begin{remark}
 We thank Ross Street for pointing us to Theorem 3.1 of \cite{johnstonealgebraic}. In that theorem, $\E$ is assumed to be an elementary topos with a natural numbers object, and the image of the natural numbers object in $\CatE$ of Corollary \ref{cor NNO one dimensional iff} is shown to be an up-to-isomorphism version of a natural numbers object in the $2$-category of toposes bounded over $\E$.
\end{remark}

\section{Subobject classifiers}\label{Section Subobject Classifiers}

We show in this section that subobject classifiers in $\E$ give rise to something similar to a subobject classifier in $\CatE$; rather than classifying monomorphisms as a subobject classifier would, the maps that are classified are monomorphisms which are also fully faithful. In this section we assume that $\E$ is lextensive and cartesian closed, so that $\CatE$ has copowers by $\mathbf{2}$ as per Theorem~\ref{thm:copower}. This means that the two-dimensional aspect of the universal property of pullbacks follows from the one-dimensional aspect, so we omit mention of it. Note $\E$ satisfying ETCS is in particular an elementary topos, and so is therefore lextensive.

\begin{definition}\label{def: full subobject classifier}
Let $\mathcal{K}$ be a $2$-category. 
\begin{enumerate}
    \item A morphism $i: A \to B$ is a \emph{full monomorphism} if for every $X\in \mathcal{K}$ the functor $\mathcal{K}\left(X, i\right): \mathcal{K}\left(X, A\right) \rightarrow \mathcal{K}\left(X, B\right)$ is fully faithful and injective on objects.
    \item Two full monomorphisms $i: A \to B$ and $i': A' \to B$ with the same codomain are said to be equivalent if there is an isomorphism $a: A \to A'$ satisfying $i' a = i$. A \emph{full subobject} of $B$ is an equivalence class of full monomorphisms into $B$.
    \item A \emph{full subobject classifier} is a full monomorphism $\underline{\top}: \underline{\mathbf{1}} \to \underline{ \Omega} $ such that for any fully faithful monomorphism $i: A \to B$, there is a unique morphism $\chi_{i}: B \to \underline{\Omega}$ making the following square a pullback.

    \begin{equation*}
        \begin{tikzcd}
            A \arrow[r, "!"] \arrow[d, "i", swap] \arrow[dr, phantom, "\lrcorner", very near start] & \underline{\mathbf{1}} \arrow[d, "\underline{\top}"] \\
            B \arrow[r, "\chi_{i}", swap] & \underline{\Omega}.
        \end{tikzcd}
    \end{equation*}
\end{enumerate}
\end{definition}

\begin{remark}\label{full subobjects generalise subobjects}
    Note that $\underline{\top}$ being a full subobject classifier is precisely to say that it is a terminal object in the category whose objects are full subobjects in $\mathcal{K}$, and whose morphisms are pullback squares. This is indeed in analogy to the universal property defining subobject classifiers, with the $2$-categorical notion of full subobjects replacing subobjects. Indeed, any monomorphism in a $1$-category $\mathcal{C}$ is fully faithful as a morphism in the discrete $2$-category on $\mathcal{C}$. As such, the notion of a full subobject classifier specialises to the notion of a subobject classifier in the setting where $\mathcal{K}$ has only identity $2$-cells. This is in contrast to other categorifications of subobject classifiers such as discrete opfibration classifiers of \cite{weber2007yoneda}. On the other hand, full-subobject classifiers in arbitrary $2$-categories are typically not subobject classifiers in their underlying categories.
   
     Note that we have not included any universal property for $2$-cells into full-subobject classifiers in Definition~\ref{def: full subobject classifier}. It is an easy exercise to check that there is a unique internal natural transformations between any parallel pair of internal functors whose codomain is an indiscrete internal category. Since the full subobject classifiers that we construct in Proposition~\ref{Cat(E) full subobject} will be indiscrete internal categories, we could have included this feature as part of the definition. We have refrained from doing so since it is not needed for Theorem~\ref{Theorem subobject classifier in E iff full subobject classifier in Cat E}, and also since doing so would lose subobject classifiers in $1$-categories as examples.
\end{remark}

\begin{prop}\label{Cat(E) full subobject}
	Suppose $\mathcal{E}$ has a subobject classifier $\top: \mathbf{1} \rightarrow \Omega$. Then $\codisc(\top): \underline{\mathbf{1}} \to \mathbf{indisc}\left(\top\right)$ is a full subobject classifier for $\CatE$.
\end{prop}

\begin{proof}
    Let $f: \X \to \Y$ be a full monomorphism. Then $f_0 : X_0 \to Y_0$ is a monomorphism in $\E$. Since $\E$ has a subobject classifier, we have a unique $\chi_{f_0}: Y_0 \to \Omega$ such that the square depicted below left is a pullback. Now, since $(-)_0 \dashv \codisc$, the adjunct of $\chi_{f_0}$ is a unique map $\chi_f: \Y \to \codisc (\Omega)$ making the square below right commute. We need to show that this square is a pullback.

    \begin{equation*}
        \begin{tikzcd}
            X_0 \arrow[r, "!"] \arrow[dr, phantom, "\lrcorner", very near start] \arrow[d, "f_0", swap]& \mathbf{1} \arrow[d, "\top"] \\
            Y_0 \arrow[r, "\exists!\chi_{f_0}", swap] & \Omega  &{}  
        \end{tikzcd}\begin{tikzcd}
            \X \arrow[r, "!"] \arrow[d, "f", swap]& \underline{\mathbf{1}} \arrow[d, "\codisc(\top)"] \\
            \Y \arrow[r, "\exists!\chi_{f}", swap] & \mathbf{indisc}(\Omega)    
        \end{tikzcd}
    \end{equation*}

 But the required square is clearly a pullback on objects, and given on morphisms as displayed below. By Proposition~\ref{arr faithful} part (2), it suffices to show that this square is a pullback. But the left square is indeed a pullback since $f: \A \to \B$ is fully faithful. The proof is complete by the pullback lemma.

\begin{equation*}
        \begin{tikzcd}[row sep = large]
            X_1 \arrow[r, "(d_0{,} d_1)"] \arrow[dr, phantom, "\lrcorner", very near start] \arrow[d, "f_1", swap]& X_0 \times X_0 \arrow[d, "f_0 \times f_0", swap] \arrow[r, "!"] \arrow[dr, phantom, "\lrcorner", very near start] & \mathbf{1} \arrow[d, "(\top {,} \top)"] \\
            Y_1 \arrow[r, "(d_0{,} d_1)", swap] & Y_0 \times Y_0 \arrow[r, "\chi_{f_0} \times \chi_{f_0}", swap]  & \Omega \times \Omega 
        \end{tikzcd}
\end{equation*}

\end{proof}

\begin{example}\label{example full subobject classifiers in Cat}
    Taking $\E = \mathbf{Set}$, the full subobject classifier in $\mathbf{Cat}$ is given by the free-living isomorphism $\mathbf{I}:= \{\bot \cong \top\}$.
\end{example}

 The proof of the converse follows easily from the adjunction $(-)_{0} \dashv \mathbf{indisc}$.

\begin{prop}\label{full soc in Cat E implies soc in E}
    Let $\E$ be a category with terminal object and pullbacks. Suppose $\CatE$ has a full subobject classifier $\underline{\top}: \underline{\mathbf{1}} \to \underline{ \Omega} $. Then $\underline{\top}_0: \mathbf{1} \to \underline{ \Omega}_0 $ is a subobject classifier for $\E$.
\end{prop}

\begin{proof}
    Let $i: A \to B$ be a monomorphism in $\E$. Then $$\codisc(i): \codisc(A) \to \codisc(B)$$ is clearly fully faithful and mono on objects; monomorphisms are closed under products and the maps $\mathbf{indisc}(X)_{1} \to \mathbf{indisc}(X)_{0} \times \mathbf{indisc}(X)_{0}$ are identities for $X \in \{A, B\}$, so that the relevant square defining fully-faithfulness is indeed a pullback. Hence, there exists a pullback square in $\CatE$ as displayed below left. Since $(-)_0$ is a right adjoint, it preserves limits and in particular pullbacks. Hence, using the fact that $(-)_0 \circ \codisc = \underline{ \mathbf{1}},$ we have the pullback square in $\E$ depicted below right. Uniqueness also follows by adjointness.

    \begin{equation*}
        \begin{tikzcd}
            \codisc(A) \arrow[r, "!"] \arrow[d, "\codisc(i)", swap] \arrow[dr, phantom, "\lrcorner", very near start] & \underline{\mathbf{1}} \arrow[d, "\underline{\top}"] \\
            \codisc(B) \arrow[r, "\phi", swap] & \underline{\Omega}&{}
        \end{tikzcd}\begin{tikzcd}
            A \arrow[r, "!"] \arrow[d, "i", swap] \arrow[dr, phantom, "\lrcorner", very near start] & \mathbf{1} \arrow[d, "\underline{\top}_0"] \\
           B \arrow[r, "\phi_0", swap] & \underline{\Omega}_0
        \end{tikzcd}
    \end{equation*}
    
\end{proof}

\begin{remark}
\label{rem proof without extensivity}
    Note that the above proof holds without the assumption that $\E$ is extensive.
\end{remark}

 We have just proven the following result.

\begin{thm}\label{Theorem subobject classifier in E iff full subobject classifier in Cat E}
    Let $\E$ be an extensive, cartesian closed category with finite limits. Then $\E$ has a subobject classifier if and only if the $2$-category $\CatE$ has a full subobject classifier. In this case, the $2$-functor $\mathbf{indisc}: \E \to \CatE$, with $\E$ being considered as a locally discrete $2$-category, preserves full-subobject classifiers. 
\end{thm}

\begin{proof}
    Combine Propositions~\ref{Cat(E) full subobject} and~\ref{full soc in Cat E implies soc in E}.
\end{proof}

\begin{remark}\label{remark booleanness and two-valuedness}
    We characterise booleanness and two-valuedness of $\E$ in terms of properties in $\CatE$. These properties follow for $\E$ from the axioms of ETCS. Booleanness is a consequence of the axiom of choice \cite{diaconescu1975axiom}, and in fact both of these properties are a consequence of well-pointedness (Proposition 7, Part VI of \cite{maclane2012sheaves}). As such, the equivalent properties that we are about to describe in $\CatE$ will also follow as a consequence of the axioms in the elementary theory of the $2$-category of small categories, which we will give in Subsection~\ref{Section ET2CSC}.
    
     Consider the two internal functors $\underline{\mathbf{1}} \to \mathbf{2}_\E$ which are the source and target of the universal $2$-cell exhibiting $\mathbf{2}_\E$ as the copower of $\underline{\mathbf{1}} \in \CatE$ by $\mathbf{2} \in \Cat$. It is easy to see that these are both full monomorphisms. Hence by Proposition~\ref{Cat(E) full subobject}, they determine internal functors $\mathbf{2}_\E \to \mathbf{indisc}(\Omega)$. At the level of objects, one of these is given by $(\top, \bot): \mathbf{1}+\mathbf{1} \to \Omega$ while the other is given by $(\bot, \top)$. Recall (Proposition 5.14 of \cite{johnstone2014topos}) that an elementary topos $\E$ is \emph{boolean} if and only if these morphisms are invertible. As such, $\E$ is boolean if and only if either (hence both) of these internal functors in $\CatE$ are codescent morphisms, since as discussed in Remark~\ref{Bourke black box} these are precisely the internal functors which are isomorphic on objects. Similarly, recall that an elementary topos is \emph{two-valued} if and only if the hom-set $\E(\mathbf{1}, \Omega)$ has exactly two morphisms, namely $\top$ and $\bot$. Hence by fully faithfulness of $\mathbf{indisc}: \E \to \CatE_{1}$, $\E$ is two-valued if and only if in $\mathbf{Cat}(\E)$ there are exactly two morphisms from the terminal object to the full subobject classifier. In this case the hom-category $\CatE(\underline{\mathbf{1}}, \mathbf{indisc}(\Omega))$ is the free-living isomorphism.
\end{remark}

\begin{remark}\label{disc Omega classifier for strict bi-sieves}
    When $\E$ is an elementary topos, the internal functor $\disc(\top): \underline{\mathbf{1}} \to \mathbf{disc}(\Omega)$ is also a classifier for a certain class of monomorphisms. These are those internal functors which are monomorphisms between objects of objects, and discrete bifibrations; a notion that can either be defined representably in $\CatE$, or internally to $\E$ by asking $f_{0}d_{k} = d_{k}f_{1}$ to be a pullback for $k \in \{0, 1\}$. We call such functors \emph{strict bi-sieves}. When $\E = \mathbf{Set}$, such functors determine a subset of the set of connected components of their codomain, and are inclusions of full subcategories on all objects in those connected components. Indeed, the proof uses the adjunction $\Pi_{0}\dashv \mathbf{disc}$ of Remark~\ref{Connected Components}. We give only a sketch of the proof, since this will not be needed for any of the results in this paper.
   
     Via $\Pi_{0}\dashv \mathbf{disc}$, a classifier $\mathbb{B} \to \mathbf{disc}(\Omega)$ corresponds to a morphism $\Pi_{0}(\B)$, which in turn corresponds to a monomorphism $f': X \to \Pi_{0}(\mathbb{B})$ in $\E$. The coequaliser diagram depicted below left is sent by $\Omega^{(-)}: \E^\text{op} \to \E$ to the equaliser diagram below right. 
    
    $$\begin{tikzcd}
        B_{1} \arrow[rr, shift left = 2, "d_{0}"]\arrow[rr, shift right =2, "d_{1}"'] && B_{0} \arrow[rr, "q_{\mathbb{B}}"]&&\Pi_{0}(\mathbb{B})&{}
    \end{tikzcd}\begin{tikzcd}
        \Omega^{B_{1}} \arrow[rr, shift right = 2, leftarrow, "\Omega^{d_{0}}"']\arrow[rr, shift left =2, leftarrow, "\Omega^{d_{1}}"] &&\Omega^{B_{0}}\arrow[rr, leftarrow, "\Omega^{q_{\mathbb{B}}}"]&& \Omega^{\Pi_{0}(\mathbb{B})}
    \end{tikzcd}$$

     But the morphisms $\Omega^{d_{k}}: \Omega^{B_{0}} \to \Omega^{B_{1}}$ for $k \in \{0, 1\}$ correspond to pullbacks of monomorphisms. As such the monomorphism $f'$ corresponds to a monomorphism $f_{0}: A_{0} \to B_{0}$ whose pullback along both $d_{0}, d_{1}: B_{1} \to B_{0}$ are the same monomorphism $f_{1}: A_{1} \to B_{1}$. These data precisely correspond to a strict bi-sieve $f: \mathbb{A} \to \mathbb{B}$. One shows that this moreover satisfies $\Pi_{0}(f) = f'$.
\end{remark}

\begin{prop}
\label{prop:extensivity given fsoc}
    Suppose $\CatE$ has finite $2$-limits, is $2$-cartesian closed and has a full subobject classifier. Then:

    \begin{enumerate}
        \item $\E$ is extensive.
        \item $\CatE$ is extensive. 
    \end{enumerate}
\end{prop}

\begin{proof}
    By Proposition~\ref{prop finite weighted limits in Cat E}, the assumptions that $\CatE$ has finite $2$-limits means that $\E$ has finite limits, and so by noting Remark~\ref{rem proof without extensivity}, we can apply Proposition~\ref{full soc in Cat E implies soc in E} and obtain a subobject classifier in $\E$. By Theorem~\ref{Theorem cartesian closedness iff}, it follows that $\E$ is cartesian closed and so $\E$ is an elementary topos and therefore extensive. By Lemma~\ref{extensive coproduct Cat(E)}, it follows that $\CatE$ is extensive.
\end{proof}

\section{The axiom of choice}\label{Section Axiom of Choice}

It is well known that the axiom of choice is equivalent to the statement that any essentially surjective on objects and fully faithful functor is part of an adjoint equivalence in $\Cat$ (\cite{freyd1990categories}, 1.364). The axiom of choice is also equivalent to the proposition that any surjective-on-objects and fully faithful functor has a section. The second of these formulations is easier to treat in the context of internal category theory. Establishing this logical equivalence is the aim of Subsection~\ref{Subsection AC in terms of internal category theory}. Subsection~\ref{AC in 2-categorical terms} will consider how the property of being epimorphic-on-objects can be expressed abstractly in the $2$-category $\mathcal{K} =\CatE$ without reference to the fact that $\mathcal{K}$ is of this form. In particular, we will show that the class of epimorphic-on-objects internal functors  in $\CatE$ is precisely the left orthogonality class with respect to the fully faithful monomorphisms. For this, we need the assumption that $\E$ has an (epi, mono)-factorisation system, which is true in any elementary (or indeed pre-)topos.

\subsection{In terms of internal category theory}\label{Subsection AC in terms of internal category theory}

\begin{definition}\label{def: External axiom of choice}
    A category $\E$ is said to satisfy the \emph{external axiom of choice} if every epimorphism $e: X \to Y$ has a section. That is, there exists a map $s: Y \to X$ satisfying $es = 1_{X}$.
\end{definition}

 We give a proof that the external axiom of choice for $\E$ is equivalent to the proposition that any epimorphic-on-objects functor that is fully faithful has a section. For this equivalence, we require that $\E$ has pullbacks and products. 

\begin{lem}\label{lemma fully faithful + epi on objects implies split}
    Let $\E$ be a category with pullbacks and products and let $e: \A \to \B$ be a fully faithful internal functor. Suppose $e_0$ has a splitting $s_0: B_0 \to A_0$. Then $s_0$ extends to an internal functor $s: \mathbb{B} \to \mathbb{A}$, with assignment on arrows given as depicted below. Moreover, $es=1_{\mathbb{B}}$.

\begin{equation*}
        \begin{tikzcd}
            B_1 \arrow[d, "(d_0{,} d_1)", swap] \arrow[dr, "s_1", dashed] \arrow[drr, "1_{B_1}", bend left] & & \\
            B_0 \times B_0  \arrow[dr, "s_0 \times s_0", swap] & A_1 \arrow[r, "e_1"] \arrow[d, "(d_0{,} d_1)", swap] \arrow[dr, phantom, "\lrcorner", very near start] & B_1 \arrow[d, "(d_0{,} d_1)"] \\
            & A_0 \times A_0 \arrow[r, "e_0 \times e_0", swap] & B_0 \times B_0
        \end{tikzcd}
\end{equation*}
       
\end{lem}

 \begin{proof}
    By construction, $s_1$ is a section of $e_1: A_1 \to B_1$ and $s: = (s_0, s_1)$ forms a morphism of the underlying graphs of $\B$ and $\A$. This morphism of graphs clearly gives a splitting of $e$. We need to prove that this is well-defined as an internal functor. We show it respects identities using the universal property of $A_1$. Compatibility with the pullback projection $e_1$ follows from the commutativity of the diagram displayed below left, while compatibility with the other pullback projection follows from the commutativity of the diagram below right, for $k \in \{0, 1\}$.

\begin{equation*}
    \begin{tikzcd}
        B_0 \arrow[r, "s_0"] \arrow[dd, "i", swap] \arrow[dr, "1_{B_0}", swap] & A_0 \arrow[r, "i"] \arrow[d, "e_0"] & A_1 \arrow[dd, "e_1"] \\
        & B_0 \arrow[dr, "i"] & \\
        B_1 \arrow[r, "s_1", swap] \arrow[rr, bend left, "1_{B_1}"] & A_1 \arrow[r, "e_1", swap] & B_1 &{}
    \end{tikzcd}\begin{tikzcd}
        B_0 \arrow[r, "s_0"] \arrow[dd, "i", swap] \arrow[dr, "1_{B_0}", swap] & A_0 \arrow[ddr, "1_{A_0}"] \arrow[r, "i"]  & A_1 \arrow[dd, "d_k"] \\
        & B_0 \arrow[dr, "s_0", swap] & \\
        B_1 \arrow[r, "s_1", swap] \arrow[ur, "d_k"]  & A_1 \arrow[r, "d_k", swap] & B_1
    \end{tikzcd}
\end{equation*}

Similarly, respect for composition also follows from the universal property of $A_{1}$ as per the calculations displayed below. This completes the proof.

\begin{equation*}
    \begin{tikzcd}
        B_2 \arrow[r, "s_2"] \arrow[ddd, "m", swap]  \arrow[dr, "\pi_k"]& A_2 \arrow[rr, "m"] \arrow[dr, "\pi_k"] & & A_1 \arrow[ddd, "d_k"] \\
        & B_1 \arrow[r, "s_1"] \arrow[d, "d_k"] & A_1 \arrow[ddr, "d_k"] & \\
        & B_0 \arrow[drr, "s_0"] && \\
        B_1 \arrow[rr, "s_1", swap] \arrow[ur, "d_k"] & & A_1 \arrow[r, "d_k", swap] & A_0 &{}
    \end{tikzcd}\begin{tikzcd}[column sep = 38, row sep = 38]
        B_2 \arrow[r, "s_2"] \arrow[dd, "m", swap] \arrow[dr, "1_{B_2}", swap] & A_2 \arrow[r, "m"] \arrow[d, "e_2"] & A_1 \arrow[dd, "e_1"] \\
        & B_2 \arrow[dr, "m"] & \\
        B_1 \arrow[r, "s_1", swap] \arrow[rr, bend left, "1_{B_1}"] & A_1 \arrow[r, "e_1", swap] & B_1
    \end{tikzcd}
\end{equation*}
\end{proof}

\begin{remark}
\label{rem: sections are lalie}
    By fully-faithfullness, $s$ can be shown to be a right adjoint equivalence right inverse to $e$. The unit $\eta: 1_\A \Rightarrow se$ is determined by $1_e$ given $e= 1_\B.e = ese$ and representable fully-faithfulness of $e: \A \to \B$. Adjointness and invertibility of $\eta$ follow from representable faithfulness and conservativity of $e$, respectively.
\end{remark}

\begin{prop}\label{TFAE with external AC epi on objects}
    Let $\E$ be a category with pullbacks. The following are equivalent:
    \begin{enumerate}
    \item The external axiom of choice holds in $\E$.
    \item Any fully faithful and epimorphism-on-objects functor internal to $\E$ has a section in the $2$-category $\CatE$.
    \end{enumerate}
\end{prop}
\begin{proof}
    Let $e: \A \to \B$ be an epi-on-objects and fully faithful functor. Assuming the external axiom of choice for $\E$, the morphism $e_{0}: A_{0} \to B_{0}$ has a splitting. The splitting for the internal functor $e: \A \to \B$ is given in Lemma~\ref{lemma fully faithful + epi on objects implies split}. 

   Conversely, assume that every epi-on-objects and fully faithful functor has a section. Let $f: X \to Y$ be an epimorphism in $\E$. The internal functor $\mathbf{indisc}(f): \mathbf{indisc}(X) \to \mathbf{indisc}(Y)$ is fully faithful and an epimorphism-on-objects and hence has a section $s: \codisc(Y) \to \codisc(X)$ giving us $s_0 : Y \to X$, a section of $f$. 
\end{proof}

\begin{example}\label{Remark trivial fibrations in Cat}
    When $\E = \s$, functors which are epi on objects and fully faithful are the right class of a weak factorisation system on $\mathbf{Cat}$, with the left class being the injective on objects functors. This factorisation system features in the canonical model structure on $\Cat$. See \cite{everaert2005model, joyal2006strong} for more on homotopical aspects of internal category theory.
\end{example}

\begin{remark}\label{Remark anafunctors}
    We briefly outline how Proposition~\ref{TFAE with external AC epi on objects} sheds light on category theory internal to categories which do not satisfy the external axiom of choice. When $\E$ does not satisfy the external axiom of choice, one often works with internal anafunctors, rather than internal functors, between internal categories so that `weak equivalences' are actually adjoint equivalences \cite{makkai1996avoiding, roberts2012internal, roberts2021elementary}. Anafunctors internal to $\E$ are typically defined in terms of covering families, an important example of which is the one generated by regular epimorphisms. In this setting, internal anafunctors $\mathbb{A} \nrightarrow \mathbb{B}$ are spans of ordinary internal functors \begin{tikzcd}
        \mathbb{A} & \mathbb{F}\arrow[l, "l"']\arrow[r, "r"]& \mathbb{B}
    \end{tikzcd} in which $l$ is fully faithful, and a regular epimorphism on objects. If regular epimorphisms are stable under pullback then internal anafunctors form the morphisms of a bicategory $\mathbf{Ana}(\CatE)$, with their composition involving pullbacks in $\CatE$. There is a canonical homomorphism of bicategories $I: \CatE \to \mathbf{Ana}(\CatE)$, which is the identity on objects and a full monomorphism between hom-categories. It views a functor as an anafunctor by taking the left leg $l: \mathbb{F} \to \mathbb{A}$ to be the identity on $\A$.
    
     If any epimorphism in $\E$ is regular and $\E$ has an (epi, mono) orthogonal factorisation system, as is the case when $\E$ is an elementary topos, then by Remark~\ref{rem: sections are lalie}, Proposition~\ref{TFAE with external AC epi on objects} says precisely that the external axiom of choice holds for $\E$ if and only if the left leg $l: \mathbb{F} \to \mathbb{A}$ in any internal anafunctor is in fact a left adjoint left inverse equivalence in $\CatE$. In this case, the homomorphism of bicategories $I: \CatE \to \mathbf{Ana}(\CatE)$ has functors between hom-categories which are essentially surjective on objects. Thus if the external axiom of choice holds for $\E$ then the $2$-category $\CatE$ is biequivalent to the bicategory $\mathbf{Ana}(\CatE)$. These observations will be generalised to appropriate $2$-categories $\mathcal{K}$ in place of $\CatE$ in Remark~\ref{Remark categorified AC implies K -> Ana(K) is a biequivalence}.
\end{remark}
\subsection{In $2$-categorical terms}\label{AC in 2-categorical terms}

 The property of being an epimorphism on objects may appear difficult to express in terms of the $2$-categorical structure of $\mathcal{K} =\CatE$, without reference to the fact that it is of this form. To fix this, we first show in Proposition~\ref{prop:liftfactorisation}, to follow, that orthogonal factorisation systems on $\E$ give rise to orthogonal factorisation systems on the $2$-category $\CatE$, as defined in \cite{day2006adjoint} and described explicitly in Remark 2.3.2 of \cite{bourke2010codescent}. This result is stated without proof in the discussion between Propositions 62 and 63 of \cite{bourke2014two}. We believe it to be of independent interest, and give a detailed proof in Appendix~\ref{appendix}. For our purposes, it will mean that epimorphism on objects internal functors can then be characterised via this left orthogonality property.

\begin{prop}
\label{prop:liftfactorisation}
    Let $(\Lagr, \R)$ be an orthogonal factorisation system on a category $\E$ with pullbacks and products. Then $$(\Lagr \text{-on-objects}, \R \text{-on-objects and fully faithful})$$ is an orthogonal factorisation system on the $2$-category $\CatE$.
\end{prop}

\begin{cor}\label{epi on objects fully mono factorisation system on Cat E}
    Let $\E$ be a category with pullbacks, products, and an orthogonal factorisation system $(\Lagr, \R)$ in which $\Lagr$ are the epimorphisms and $\R$ are the monomorphisms. Then $(\Lagr', \R')$ is an orthogonal factorisation system on $\CatE$, where $\Lagr'$ is the class of internal functors which are epi-on-objects, and $\R'$ is the class of full monomorphisms.
\end{cor}

\begin{proof}
    By Proposition~\ref{prop:liftfactorisation} there is an orthogonal factorisation system $(\Lagr', \R')$ on the $2$-category $\CatE$ in which $\Lagr'$ is as required and $\R'$ is the class of internal functors which are both fully faithful and given by monomorphisms on objects. But as discussed in the beginning of Remark~\ref{foreshadowing remark full subobjects}, such internal functors are precisely the full monomorphisms in $\CatE$. 
\end{proof}

\begin{remark}\label{Remark so/full mono lit review}
    The factorisation system on $\CatE$ obtained in Corollary~\ref{epi on objects fully mono factorisation system on Cat E} is an internal version of the factorisation system constructed in $\Cat$ via kernels and quotients, in 5.2 of \cite{bourke2014two}. The class of full monomorphisms, and its left orthogonality class, are respectively called \emph{chronic} and \emph{acute} in 1.1 and 1.4 of \cite{street1982two}. Although pullback stability of maps in $\Lagr\subseteq \E$ is not needed in the proof of Proposition~\ref{prop:liftfactorisation}, if this class is pullback stable then so is the class $\Lagr' \subseteq \CatE$. In particular, if $\E$ is a regular category then $\CatE$ is a regular $2$-category in the sense of 1.19 in \cite{street1982two}.
\end{remark}

\begin{definition}\label{Definition acute}
    (\cite{street1982two}) A morphism in a $2$-category $\mathcal{K}$ which is left orthogonal to all fully faithful monomorphisms in $\mathcal{K}$ will be called \emph{acute}.
\end{definition}

\begin{definition}\label{definition categorified axiom of choice}
    Say that a $2$-category $\mathcal{K}$ satisfies the \emph{categorified axiom of choice} if any acute fully faithful morphism has a section.
\end{definition}

 Putting these results together gives the following reformulation of the external axiom of choice in $\E$ in terms of the $2$-categorical structure of $\mathbf{Cat}(\E)$.

\begin{thm}\label{Theorem external AC iff every ff + left orhtogonal to full mono is a left adjoint left inverse equivalence}
     Let $\E$ be a category with pullbacks, products and an (epi, mono)-orthogonal factorisation system. Then the following are equivalent.
    \begin{enumerate}
    \item The category $\E$ satisfies the external axiom of choice.
    \item The $2$-category $\CatE$ satisfies the categorified axiom of choice.
    \end{enumerate}
\end{thm}

\begin{proof}
    Proposition~\ref{TFAE with external AC epi on objects} established the logical equivalence between the external axiom of choice in $\E$ and an analogue of the categorified axiom of choice for $\CatE$ with `epi-on-objects' in place of acute. But Corollary~\ref{epi on objects fully mono factorisation system on Cat E} ensures that being an epimorphism on objects characterises acute morphisms in $\CatE$.
\end{proof}

\begin{remark}\label{Remark categorified AC implies K -> Ana(K) is a biequivalence}
    The discussion in Remark~\ref{Remark anafunctors} is also possible to rephrase in $2$-categorical terms, rather than in terms of internal category theory. Let $\mathcal{K}$ be a $2$-category with pullbacks and suppose that acute morphisms are stable under pullback in $\mathcal{K}$. Define an \emph{anamorphism} in $\mathcal{K}$ to be a span whose left leg is acute and fully faithful. Then there is a bicategory $\mathbf{Ana}(\mathcal{K})$ defined in the usual way. There is also a homomorphism of bicategories $I: \mathcal{K} \to \mathbf{Ana}(\mathcal{K})$ which is given by the identity on objects and full monomorphisms between hom-categories. If the categorified axiom of choice holds in $\mathcal{K}$, then $I$ moreover has functors between hom-categories which are essentially surjective on objects. Hence in this case $I$ a biequivalence, exhibiting morphism composition as a strictification of anamorphism composition.
\end{remark}

\begin{remark}

  We thank Richard Garner for observing that when $\E$ is regular, acuteness of fully faithful internal functors is equivalent to the simpler property of being a regular epimorphism. It is clear that if $\E$ has products, then since $(-)_{0}: \CatE_{1} \to \E$ is a left adjoint it preserves regular epimorphisms. Conversely, if $f_{0}: \A_{0} \to \B_{0}$ is a regular epimorphism and $\E$ is a regular category then $f_{0}$ is the coequaliser of its kernel pair in $\E$. Then $f_{1}: A_{1} \to B_{1}$ is also a regular epimorphism, since $f: \A \to \B$ is fully faithful and regular epimorphisms are closed under products and stable under pullback in $\E$. One verifies that $f$ is the coequaliser of its kernel pair in $\CatE$ using the universal property of the coequalisers $f_{0}$ and $f_{1}$ in $\E$; we leave these details to the interested reader.
\end{remark}

\section{Comparing ETCS to ET2CSC}\label{Section comparing ETCS to ET2CSC}

 We collect the main results of previous sections and characterise $2$-categories of the form $\CatE$ when $\E$ is a model of the elementary theory of the category of sets. Our characterisation of such $2$-categories is in $2$-categorical terms, rather than in terms of category theory internal to the discrete objects of $\mathcal{K}$. The theory of such $2$-categories is again elementary, although we refrain from providing an explicit first order presentation as is done for ETCS on \cite{nlab:fully_formal_etcs}. Following this, in Subsection~\ref{Subsection Functorial ET2CSC} we describe relationships between different models of ET2CSC, and establish a `Morita biequivalence' between ETCS and ET2CSC.

\subsection{A characterisation of $\CatE$ when $\E$ is a model of ETCS}\label{Section ET2CSC}

\begin{definition}
\label{def:ET2CSC}
    We say that the $2$-category $\mathcal{K}$ models the \emph{elementary theory of the $2$-category of small categories} (ET2CSC) if the following properties hold:

    \begin{enumerate}
        \item It satisfies the conditions listed in Proposition~\ref{Bourke characterisation of Cat E}
        \item It has a terminal object.
        \item It is cartesian closed.
        \item It is $2$-well-pointed, in the sense of Definition~\ref{def:2-well-pointed}.
        \item It has a natural numbers object, in the sense of Definition~\ref{definition NNO} part (2).
        \item It has a full subobject classifier, in the sense of Definition~\ref{def: full subobject classifier} part (3).
        \item It satisfies the categorified axiom of choice, in the sense of Definition~\ref{definition categorified axiom of choice}.
    \end{enumerate}
\end{definition}

 We are now ready to combine the results so far and prove our first main result.

\begin{thm}\label{thm:ET2CSC}

\begin{enumerate}
    \item Let $\E$ be a category. Then $\E$ models the elementary theory of the category of sets if and only if $\CatE$ models the elementary theory of the $2$-category of small categories, and in this case $\E \simeq \mathbf{Disc}(\CatE)$.
    \item Conversely, let $\mathcal{K}$ be a $2$-category. Then $\mathcal{K}$ models the elementary theory of the $2$-category of small categories if and only if $\mathbf{Disc}\left(\mathcal{K}\right)$ models the elementary theory of the category of sets, and in this case $\mathcal{K} \simeq \mathbf{Cat}(\mathbf{Disc}\left(\mathcal{K}\right))$.
\end{enumerate} 
\end{thm}

\begin{proof}
    Proposition~\ref{Bourke characterisation of Cat E} gives the correspondence between pullbacks in $\E$ and the first item of Definition~\ref{def:ET2CSC}, as well as the equivalences $\E \simeq \mathbf{Disc}(\CatE)$ and $\mathcal{K} \simeq \mathbf{Cat}(\mathbf{Disc}\left(\mathcal{K}\right))$. We describe how the results in this paper so far give correspondences between the various other properties of ETCS and ET2CSC.
   
     The correspondence for terminal objects is given in Proposition~\ref{prop finite weighted limits in Cat E}, and the correspondence for cartesian closedness is Theorem~\ref{Theorem cartesian closedness iff}. Herein, assume that the category $\E$ (resp. the $2$-category $\mathcal{K}$) satisfies the properties mentioned so far. Assuming additionally that $\E$ has a subobject classifier makes $\E$ into an elementary topos; in particular it is extensive and so by Theorem~\ref{Theorem subobject classifier in E iff full subobject classifier in Cat E}, $\CatE$ has a full subobject classifier. Conversely, assuming that $\mathcal{K} \simeq \Cat(\mathbf{Disc}(\mathcal{K}))$ has a full subobject classifier means that $\E$ is extensive by Proposition~\ref{prop:extensivity given fsoc} and so by Theorem~\ref{Theorem subobject classifier in E iff full subobject classifier in Cat E}, we get the other direction of this correspondence. Note also that by Theorem~\ref{thm:copower}, under these assumptions $\mathcal{K}$ (resp. $\CatE$) has copowers by $\mathbf{2}$. 
     
     The correspondence between well-pointedness and $2$-well-pointedness is Theorem~\ref{Theorem generator in E vs generator in CatE}. The correspondence for natural numbers objects is Theorem~\ref{thm NNO iff}. Finally, the correspondence between the axiom of choice and the categorified axiom of choice is Theorem~\ref{Theorem external AC iff every ff + left orhtogonal to full mono is a left adjoint left inverse equivalence}. This last correspondence uses the epi-mono factorisation system on $\E$ (resp. $\mathbf{Disc}(\mathcal{K})$), which exists since by this stage this category is an elementary topos.
\end{proof}

 Theorem~\ref{thm:ET2CSC} will be built upon further in Subsection~\ref{Subsection Functorial ET2CSC}, where we will define $2$-categories whose objects are models of ETCS and ET2CSC respectively, and prove that these two $2$-categories are biequivalent in Theorem~\ref{Theorem morita equivalence between ETCS and ET2CSC}.

\begin{remark}\label{Remark 2D aspects for ET2CSC follow from 1D aspects}
    Assuming that $\mathcal{K}$ satisfies the conditions listed in Proposition~\ref{Bourke characterisation of Cat E}, the one-dimensional aspects of the remaining conditions in ET2CSC are enough to imply that $\mathbf{Disc}(\mathcal{K})$ satisfies ETCS, and hence that $\mathcal{K}$ satisfies the two-dimensional aspects of ET2CSC. In particular, the theory can be simplified by removing the two-dimensional aspect of cartesian closedness, the faithfulness on $2$-cells aspect of $2$-well-pointedness, the two-dimensional aspect of the universal property of natural numbers objects, and the two-dimensional aspect of left orthogonality in the definition of acute maps. Indeed, as discussed in Remark~\ref{full subobjects generalise subobjects}, we could have also included a two-dimensional universal property in our definition of a full-subobject classifier. Such a definition would demand a representing object $\underline{\Omega}$ for the $2$-functor $\mathcal{K}^\text{op} \to \mathbf{Cat}$ which sends an object $X$ to the indiscrete category on the set of full subobjects into $X$, and acts on morphisms via pullback. We chose not to give such a definition so that we retained ordinary subobject classifiers as examples.
\end{remark}

\subsection{Morphisms of models of ET2CSC}\label{Subsection Functorial ET2CSC}

 The notion of what a morphism of models of ETCS or of ET2CSC should be is clear from the description of these theories, but we spell it out in detail in Definition~\ref{Definition conditions for a et2csc morphism}, to follow. The aim of this Subsection is to extend Theorem~\ref{thm:ET2CSC} to a correspondence between morphisms of models of the two theories, and to show that they have biequivalent $2$-categories of models.

\begin{definition}\label{Definition conditions for a et2csc morphism}

    \begin{enumerate}
        \item \label{morphism of ETCS} Let $\E$ and $\E'$ be categories modelling ETCS. An \emph{ETCS-morphism} is a functor $F: \E \to \E'$ which preserves finite limits, internal homs, the subobject classifier, and the natural numbers object.

        \item \label{morphism of ET2CSC} Let $\mathcal{K}$ and $\mathcal{K}'$ be $2$-categories modelling ET2CSC. An \emph{ET2CSC-morphism} is a $2$-functor $\mathbf{F}: \mathcal{K} \to \mathcal{K}'$ which preserves pullbacks, powers by $\mathbf{2}$, codescent objects of cateads, the terminal object, internal homs, the full-subobject classifier, and the natural numbers object.
    \end{enumerate}
\end{definition}

\begin{prop}\label{Bourke characterisation of Cat F}
    (Theorem 4.28 of \cite{bourke2010codescent}) Let $F: \E \rightarrow \E'$ be a pullback preserving functor. Then $\mathbf{Cat}\left(F\right)$ preserves pullbacks, powers by $\mathbf{2}$ and codescent objects of cateads and there is a natural isomorphism $F \cong \mathbf{Disc}\circ \Cat (F)$. Conversely, if a $2$-functor $G: \mathbf{Cat}\left(\E\right) \rightarrow \mathbf{Cat}\left(\E'\right)$ preserves pullbacks, powers by $\mathbf{2}$ and codescent objects of cateads, then $\mathbf{Disc}\left(G\right): \mathbf{Disc}\circ \mathbf{Cat}\left(\E\right) \rightarrow \mathbf{Disc}\circ \mathbf{Cat}\left(\E'\right)$ preserves pullbacks and there is a $2$-natural isomorphism $G \cong \mathbf{Cat}\circ \mathbf{Disc}\left(G\right)$.
\end{prop}

\begin{remark}
    By Proposition~\ref{Bourke characterisation of Cat F} an ET2CSC-morphism is isomorphic to one of the form $\mathbf{F} = \Cat(F)$ for some pullback preserving functor $F: \E \to \E'$. As such, we will continue this section assuming that $\mathbf{F} \cong \Cat(F)$ for some such $F: \E \to \E'$. Note that since ET2CSC-morphisms preserve pullbacks, terminal objects and powers by $\mathbf{2}$, they preserves all $2$-limits. The reason that well-pointedness, the axiom of choice and their respective analogues do not feature in Definition~\ref{Definition conditions for a et2csc morphism} is that these are properties rather than structure to be preserved. In any case, logical functors preserve epimorphisms and the terminal object, and once we show that $\mathbf{Disc}(F)$ for a morphism of models of ET2CSC is a logical functor, it will follow in Corollary~\ref{corollary et2csc morphisms preserve coproducts and copowers by 2} that $F$ also preserves coproducts, copowers by $\mathbf{2}$, and acute morphisms.
\end{remark}

\begin{thm}\label{theorem characterisation of ET2CSC morphisms}
    A $2$-functor $\mathbf{F}: \mathcal{K} \to \mathcal{K}'$ between categories satisfying ET2CSC is an ET2CSC-morphism if and only if is is of the form $\mathbf{F} \cong \Cat(F)$ for some $F:\E \to \E'$ where $F$ is an ETCS-morphism. 
\end{thm}

 We prove this through a series of lemmata. In these, we repeatedly use the fact that $(-)_0$ is a $2$-natural transformation from the $2$-functor $\Cat(-): \mathbf{Lex} \to \mathbf{Lex}$ to the identity on $\mathbf{Lex}$. Here $\mathbf{Lex}$ denotes the $2$-category whose objects are categories with finite limits, whose morphisms are functors that preserve finite limits, and whose $2$-cells are arbitrary $2$-natural transformations. Similarly, we use that $\disc: 1_\mathbf{Lex} \to \Cat(-)$ is a $2$-natural transformation and that $\mathbf{indisc}: 1_\mathbf{Lex} \to \Cat(-)$ is a pseudonatural transformation. See \cite{miranda2022internal} for proofs of these properties, although we will address preservation of the terminal object in Lemma~\ref{Lemma Cat F preserves terminal} for completeness. Throughout these proofs, suppose that our $2$-functor $\mathbf{F}$ preserves pullbacks, powers by $\mathbf{2}$ and codescent objects of cateads, so that it is of the form $\mathbf{F} = \Cat(F)$ for some $F: \E \to \E'$.

\begin{lem}\label{Lemma Cat F preserves terminal}
    A pullback preserving functor $F: \E \to \E'$ preserves the terminal object if and only if  $\Cat(F): \Cat(\E) \to \Cat(\E')$ preserves the terminal object. 
\end{lem}

\begin{proof}

    Suppose that for any $A,B \in \E$, we have $F\mathbf{1} \cong \mathbf{1}'.$ Then by $2$-naturality of $\disc$

    \begin{equation*}
        \Cat(F)(\underline{\mathbf{1}})= \Cat(F)(\disc(\mathbf{1})) = \disc(F(\mathbf{1})) \cong \disc(\mathbf{1}') = \underline{\mathbf{1}}'. 
    \end{equation*}

     Conversely, suppose $\Cat(F)(\underline{\mathbf{1}}) \cong \underline{\mathbf{1}}'.$ By $2$-naturality of $(-)_0$, we have 

    \begin{equation*}
        F\mathbf{1}  = F(\underline{\mathbf{1}})_0
        =(\Cat(F)\underline{\mathbf{1}})_0 \cong (\underline{\mathbf{1}}')_{0}
         = \mathbf{1}'.
    \end{equation*}
\end{proof}

 In Lemma~\ref{lemma cat f preserves exponentials}, to follow, we denote exponentials in $\E$ as $[X, Y]$ rather than $Y^X$, for ease of readability. Similarly, we denote exponentials in $\CatE$ as $\underline{\pmb{[}\X,\Y\pmb{]}}$.

\begin{lem}\label{lemma cat f preserves exponentials}
    Suppose $F: \E \to \E'$ preserves finite limits. Then $F[A,B]\cong[FA, FB]'$ for all $A, B \in \E$ if and only if $\Cat(F)\underline{\pmb{[}\mathbb{X}, \mathbb{Y}\pmb{]}} \cong \underline{\pmb{[}\Cat(F)\mathbb{X}, \Cat(F) \mathbb{Y}\pmb{]}}'$ for all $\mathbb{X}, \mathbb{Y} \in \CatE$. 
\end{lem}

\begin{proof}

    Suppose $F[A,B]\cong[FA, FB]'$ and recall that exponentials in $\CatE_{1}$ are constructed in $[\Delta_{\leq 3}^\text{op}, \E]$. The proof that $\Cat(F)\underline{\pmb{[}\mathbb{X}, \mathbb{Y}\pmb{]}} \cong \underline{\pmb{[}\Cat(F)\mathbb{X}, \Cat(F) \mathbb{Y}\pmb{]}}'$ follows from the chain of isomorphisms in $[\Delta_{\leq 3}^\text{op}, \E]$ depicted below.

    \begin{align*}
        \Cat(F)\underline{\pmb{[}\X,\Y\pmb{]}}(-) & = F \int_{[n] \in \Delta_{\le 3}} \prod\limits_{\phi \in \Delta(-, n)}[X_n,Y_n] & \text{ definition of exponentials in $\CatE$,}\\
        & \cong \int_{[n] \in \Delta_{\le 3}} F\prod\limits_{\phi \in \Delta(-, n)}[X_n,Y_n] & \text{ the end is a finite limit,}\\
        & \cong \int_{[n] \in \Delta_{\le 3}} \prod\limits_{\phi \in \Delta(-, n)}F[X_n,Y_n] & \text{ each hom of $\Delta_{\leq 3}$ is finite,} \\
        & \cong \int_{[n] \in \Delta_{\le 3}} \prod\limits_{\phi \in \Delta(-, n)}[FX_n,FY_n]' & \text{ $F$ preserves exponentials,} \\
        & = \underline{\pmb{[}\Cat(F)\X, \Cat(F)\Y\pmb{]}}'(-) & \text{ by definition of $\Cat(F).$} 
    \end{align*}

     Conversely, suppose that for any $\X, \Y \in \CatE$, we have $$\Cat(F)\underline{\pmb{[}\mathbb{X}, \mathbb{Y}\pmb{]}} \cong \underline{\pmb{[}\Cat(F)\mathbb{X}, \Cat(F) \mathbb{Y}\pmb{]}'}.$$ In Theorem~\ref{Theorem cartesian closedness iff}, we showed that $[A,B] = ([\disc(A), \disc(B)])_0$. Let $A, B \in \E$. Then 

    \begin{align*}
        F[A,B] &= F\underline{\pmb{[}\disc(A), \disc(B)\pmb{]}}_0 \\
        & = (\Cat(F)\underline{\pmb{[}\disc(A), \disc(B)\pmb{]}})_0 \\
        & \cong \underline{\pmb{[}\Cat(F)\disc(A), \Cat(F) \disc(B)\pmb{]}}'_0 \\
        & = \underline{\pmb{[}\disc(FA), \disc(FB)\pmb{]}}'_0 \\
        & = [FA, FB]'
    \end{align*}

\end{proof}

\begin{lem}\label{Cat F preserves full subobject classifiers}
    $F \Omega \cong \Omega'$ if and only if $\Cat(F)( \underline{\Omega}) \cong \underline{\Omega}'.$ 
\end{lem}

\begin{proof}
    In Section~\ref{Section Subobject Classifiers}, we characterised the full subobject classifier of $\Cat(\E)$ in terms of the subobject classifier in $\E$, with the full subobject classifier being given by $\underline{\Omega} := \codisc(\Omega).$ 

   Assume that $F: \E \to \E'$ preserves the subobject classifier. Then there is the following chain of isomorphisms in $\Cat(\E')$, with the first being given by pseudonaturality of $\mathbf{indisc}$ in $F$ and the second being given by the isomorphism up to which $F$ preserves the subobject classifier.

    \begin{equation*}
        \Cat(F)(\underline{\Omega}) = \Cat(F)(\codisc(\Omega)) \cong \codisc(F \Omega ) \cong \codisc(\Omega') = \underline{\Omega}'.
    \end{equation*}

     Conversely, suppose that $\Cat(F)\underline{\Omega} \cong \underline{\Omega}'.$ Then the calculation below demonstrates that $F: \E \to E'$ also preserves the subobject classifier.

    \begin{equation*}
        F\Omega = F(\codisc(\Omega))_0 = F (\underline{\Omega})_0 = (\Cat(F)\underline{\Omega})_0 \cong (\underline{\Omega}')_0 = (\codisc(\Omega'))_0 = \Omega'.
    \end{equation*}
\end{proof}

\begin{cor}\label{corollary et2csc morphisms preserve coproducts and copowers by 2}
    If $\mathbf{Cat}\left(F\right)$ preserves pullbacks, powers by $\mathbf{2}$ and codescent objects of cateads, then $\mathbf{Cat}\left(F\right)$ preserves coproducts, copowers by $\mathbf{2}$, and acute morphisms as in Definition~\ref{Definition acute}.
\end{cor}

\begin{proof}
    By Proposition~\ref{Bourke characterisation of Cat F}, Lemma~\ref{Lemma Cat F preserves terminal}, Lemma~\ref{lemma cat f preserves exponentials} and Lemma~\ref{Cat F preserves full subobject classifiers}, it follows that $F: \E \to \E'$ is a logical functor. But logical functors preserve coproducts (Corollary 2.2.10 part (i) A2.2 \cite{johnstone2002sketches}), and coproducts in $\mathbf{Cat}\left(\E\right)$ are computed in $[\Delta^\text{op}, \E]$ so $\mathbf{Cat}\left(F\right)$ also preserves coproducts. Similarly, $\mathbf{Cat}\left(F\right)$ preserves copowers by $\mathbf{2}$ since these are built in $\E$ out of coproducts, terminal objects and products, all of which $F$ preserves. Finally, by Corollary~\ref{epi on objects fully mono factorisation system on Cat E}, acute morphisms in $\CatE$ are precisely the epimorphism-on-objects internal functors, and logical functors also preserve epimorphisms.
\end{proof}

\begin{lem}\label{preservation of NNO}
    $F(N) \cong N'$ if and only if $\Cat(F)(\underline{N}) \cong \underline{N}'$.
\end{lem}

\begin{proof}
    Similar to the proof of Lemma \ref{Cat F preserves full subobject classifiers}, but with $\disc$ in place of $\codisc$.
\end{proof}

 We now describe how these results combine to prove Theorem~\ref{theorem characterisation of ET2CSC morphisms}.

\begin{proof}
    (Theorem~\ref{theorem characterisation of ET2CSC morphisms}).
   
     The correspondence between preservation of pullbacks in $\E$, and preservation of pullbacks, powers by $\mathbf{2}$, and codescent objects of cateads in $\mathcal{K}$ is part of Bourke's result recalled in Proposition~\ref{Bourke characterisation of Cat F}. The correspondence for preservation of terminal objects is shown in Lemma~\ref{Lemma Cat F preserves terminal}, while the correspondence for preservation of exponentials is shown in Lemma~\ref{lemma cat f preserves exponentials}. The correspondence between preservation of subobject classifiers in $\E$ and full subobject classifiers in $\mathcal{K}$ is shown in Lemma~\ref{Cat F preserves full subobject classifiers}. Finally, the correspondence between preservation of natural numbers objects is shown in Lemma~\ref{preservation of NNO}. 
\end{proof}

 Theorem~\ref{Theorem morita equivalence between ETCS and ET2CSC}, to follow, says that ETCS and ET2CSC have biequivalent $2$-categories of models. This is the sense in which we claim to have categorified ETCS, and provided a foundation of mathematics that captures the structural aspects of categories. In contrast, ETCS is a foundation which axiomatises the structural properties of sets.

\begin{thm}\label{Theorem morita equivalence between ETCS and ET2CSC}
    Let $\mathbf{ETCS}$ denote the $2$-category whose objects are categories modelling ETCS, whose morphisms are ETCS morphisms, and whose $2$-cells are natural isomorphisms. Let $\mathbf{ET2CSC}$ denote the $2$-category whose objects are $2$-categories modelling ET2CSC, whose morphisms are ET2CSC morphisms, and whose $2$-cells are $2$-natural isomorphisms. Then there is a biequivalence as depicted below.

    $$\begin{tikzcd}
        \mathbf{ETCS} \arrow[rr,shift right =2, "\mathbf{Cat}(-)"']&\sim&\mathbf{ET2CSC}\arrow[ll, shift right = 2, "\mathbf{Disc}(-)"']
    \end{tikzcd}$$
\end{thm}

\begin{proof}
    The required biequivalence is a restriction of the one in Theorem 4.28 of \cite{bourke2010codescent}. The fact that it restricts as required follows from Theorem~\ref{thm:ET2CSC} and Theorem~\ref{theorem characterisation of ET2CSC morphisms}.
\end{proof}

\section{Conclusions and future directions}\label{Section conclusion}

 In this paper we have extended Bourke's characterisation of $2$-categories of the form $\CatE$ of internal categories, functors and natural transformations for $\E$ a category with pullbacks (Proposition~\ref{Bourke characterisation of Cat E}), and his characterisation of $2$-functors of the form $\Cat (F)$ for pullback preserving functors $F: \E \to \E'$ (Proposition~\ref{Bourke characterisation of Cat F}). Specifically, we have characterised $2$-categories of the same form $\CatE$, where $\E$ now models Lawvere's elementary theory of the category of sets (Theorem~\ref{thm:ET2CSC}), and we have also characterised $2$-functors of the form $\Cat (F)$ where $F: \E \to \E'$ preserves the structure in ETCS (Theorem~\ref{theorem characterisation of ET2CSC morphisms}). In particular, we have done so in a way that such $2$-categories $\CatE$ can be finitely axiomatised in first order logic, without presupposing an ambient set theory. For these reasons we have called the theory of such $2$-categories `the elementary theory of the $2$-category of small categories', or ET2CSC. These results build upon Bourke's work to show that ETCS and ET2CSC have biequivalent $2$-categories of models (Theorem~\ref{Theorem morita equivalence between ETCS and ET2CSC}). To the extent that ETCS provides a structural foundation by axiomatising the category structure of sets and functions, ET2CSC provides a structural foundation by axiomatising the $2$-category structure of categories, functors and natural transformations.

ET2CSC also has the feature that it can be expressed in purely $2$-categorical terms, without reference to the fact that its models are of the form $\CatE$, up to equivalence. An important step towards this is Corollary~\ref{cor:generator}, in which we show that generating families in lextensive $\E$ give rise to generating families in $\CatE$. This motivated the notion of $2$-well-pointedness, introduced in Definition~\ref{def:2-well-pointed} part (2), which is a key ingredient in ET2CSC. Another key ingredient in this axiomatisation is the concept of a `full subobject', which is an abstraction of functors which include full subcategories determined by some subset of the objects of their codomain. Classifiers for full-subobjects were introduced in Definition~\ref{def: full subobject classifier}, and such classifiers in $\CatE$ were shown in Theorem~\ref{Theorem subobject classifier in E iff full subobject classifier in Cat E} to be tantamount to subobject classifiers in $\E$. Meanwhile, maps which are left orthogonal to these full subobjects, the so called acute maps of \cite{street1982two}, played a role in expressing the categorified axiom of choice abstractly rather than in terms of internal category theory, in Theorem~\ref{Theorem external AC iff every ff + left orhtogonal to full mono is a left adjoint left inverse equivalence}. The correspondence between specific properties of $\E$ and analogous properties of $\CatE$ often requires much less to be assumed than the remaining properties in ETCS, and we have presented our proofs accordingly so that the various intermediate results may be applied in greater levels of generality. In particular, we think the intermediate results Corollary~\ref{cor:generator} and Proposition~\ref{prop:liftfactorisation} may be of independent interest.

 It is also of interest to establish sufficient elementary conditions on $\E$ for finite $2$-colimits to exist in $\CatE$. Coproducts and copowers by $\mathbf{2}$ were treated in section~\ref{Section Generators}, while Remark~\ref{remark coequalisers need NNOs} recorded that $\E$ being an elementary topos is insufficient for $\CatE$ to have coequalisers, or even coinserters. We conjecture that being an elementary topos with a natural numbers object \emph{is} sufficient for coequalisers, and hence finite $2$-colimits, to exist in $\CatE$. Indeed, coequalisers of a parallel pair of functors depicted below left are constructed in $\Cat$ using not just coequalisers in $\s$, by also lists of morphisms in $\mathcal{D}$. Specifically, consecutive morphisms $(g_{n+1}, g_{n})$ in such lists are not already composable in $\mathcal{D}$, but rather their intermediate objects $d_{1}(g_{n+1})$ and $d_{0}(g_{n})$ must be identified by the coequaliser in $\s$ depicted below right. Coequalisers of a general parallel pair of functors are constructed similarly, but moreover involve generating a congruence from $Ff\sim Gf$ for morphisms $f$ in their domain, and then quotienting by this congruence.

$$\begin{tikzcd}
    \mathbf{disc}(X)\arrow[rr, shift left = 2, "F"]\arrow[rr, shift right = 2, "G"']&& \mathcal{D}&{}
\end{tikzcd} \begin{tikzcd}
    X\arrow[rr, shift left = 2, "F_{0}"]\arrow[rr, shift right = 2, "G_{0}"']&& \mathcal{D}_{0}
\end{tikzcd}$$ 

In future work \cite{HughesMiranda20242categories} we extend the theory developed here to incorporate the axiom of replacement in the framework of a $2$-category of categories. Another interesting direction for future research would be to reformulate other set theoretical conditions such as the continuum hypothesis, or large cardinal axioms, in terms of the $2$-categorical structure of $\CatE$. On the other hand, a related but different direction for future research could be axiomatising $2$-categories of the form $\CatE$ when $\E$ satisfies Giraud's axioms for Grothendieck toposes (Proposition 6.1.01 of \cite{lurie2009higher}), Frey's axioms for realisability toposes \cite{frey2019characterizing}, or Kock's axioms for smooth toposes \cite{kock2006synthetic}. Finally, one could try to extend our work to higher categorical settings by axiomatising the three dimensional structure that small double categories and double functors underlie \cite{bohm2020gray}. Similarly, higher dimensional structures comprising Kan complexes or quasicategories are already an active area of research \cite{riehl2022elements, stenzel2024infty2category}. 

\appendix
\section{Proof of Proposition~\ref{prop:liftfactorisation}}
\label{appendix}

\begin{notation}\label{notation L' R' ofs on Cat E}

    \begin{enumerate}
    \item For $s, f$ morphisms in a $2$-category $\mathcal{K}$, write $s \perp_{1} f$ if for any commutative square in $\mathcal{K}$ as depicted below left, there is a unique morphism $u: B \to X$ satisfying $fu = q$ and $us = p$.
    \item Write $s \perp f$ if $s \perp_{1} f$ and moreover for any commutative pair of $2$-cells as depicted below right, with $u: B \to X$ and $u': B \to X$ the corresponding morphisms induced by $fp=qs$ and $fp'=q's$ respectively, there is a unique $2$-cell $\gamma: u \Rightarrow u'$ satisfying $f.\gamma = \beta$ and $\gamma.s = \alpha$. 

    \begin{equation*}
        \begin{tikzcd}[row sep = 30, column sep = 30]
            A \arrow[r, "p"] \arrow[d, "s", swap] & X \arrow[d, "f"] \\
            B \arrow[r, "q", swap] & Y&{}
        \end{tikzcd}\begin{tikzcd}[row sep = 30, column sep = 30]
            A \arrow[r, bend left = 30,"p"name=A]\arrow[r, bend right = 30,"p'"'name=B] \arrow[d, "s", swap]
            & X \arrow[d, "f"] \\
            B \arrow[r, bend left = 30, "q" name=C]\arrow[r, bend right = 30, "q'"' name=D] & Y
            \arrow[from=A, to=B, Rightarrow, "\alpha"]
            \arrow[from=C, to=D, Rightarrow, "\beta"]
        \end{tikzcd}
    \end{equation*}
    \item For $\mathcal{M}_{1}$, $\mathcal{M}_{2}$ classes of morphisms in a category $\mathcal{C}$, write $\mathcal{M}_{1} \perp_{1} \mathcal{M}_{2}$ if for every $s \in \mathcal{M}_{1}$ and $f \in \mathcal{M}_{2}$, we have $s \perp_{1} f$.
    \item For $\mathcal{M}_{1}$, $\mathcal{M}_{2}$ classes of morphisms in a $2$-category $\mathcal{K}$, write $\mathcal{M}_{1} \perp \mathcal{M}_{2}$ if for every $s \in \mathcal{M}_{1}$ and $f \in \mathcal{M}_{2}$, we have $s \perp f$.
    \item For $(\Lagr, \R)$ an orthogonal factorisation system on $\E$, let $\Lagr'$ denote the class of internal functors which are $\Lagr$-on-objects, and let $\R'$ denote the class of internal functors which are $\R$-on-objects and fully faithful. 
    \end{enumerate}
\end{notation}

 It is clear that both $\Lagr'$ and $\R'$ contain all isomorphisms of internal categories and are closed under composition, since these properties hold for the classes of morphisms $\Lagr$ and $\R$ in $\E$, and for the class of fully faithful functors in $\CatE$. By Lemma 2.2 of \cite{bousfield1977constructions}, it therefore suffices to show that the following properties hold to establish that $(\Lagr', \R')$ is an orthogonal factorisation system on the category $\CatE_{1}$.

\begin{itemize}
    \item $\Lagr' \perp_{1} \R'$. 
    \item Any internal functor $f: \mathbb{X} \to \mathbb{Y}$ admits a factorisation $f=rl$ with $l \in \Lagr'$ and $r \in \R'$. 
\end{itemize} 

 If moreover $\Lagr' \perp \R'$, then $(\Lagr', \R')$ is an orthogonal factorisation system on the $2$-category $\CatE$. We prove $\Lagr' \perp \R'$ in Lemma~\ref{lemma orthogonality}, and the existence of an appropriate factorisation in Lemma~\ref{lemma factorisation}.

\begin{lem}\label{lemma orthogonality}
   $\Lagr' \perp \R'$.
\end{lem}

\begin{proof}
    We first prove the one-dimensional aspect of orthogonality. Consider a diagram in $\CatE$ as depicted below, in which $s \in \Lagr'$ and $f \in \R'$.

    $$\begin{tikzcd}
            \A \arrow[r, "p"] \arrow[d, "s", swap] & \X \arrow[d, "f"] \\
            \B \arrow[r, "q", swap] & \Y
    \end{tikzcd}$$

     Apply the functor $(-)_0: \CatE_{1} \to \E$ to get a commutative square as depicted below left in $\E$, in which the unique lift exists as $s_0 \in \Lagr$ and $f_0 \in \R$. We define $u_1: B_1 \to X_1$ by the universal property of $X_1$, as depicted below right.

$$\begin{tikzcd}
            A_0 \arrow[rr, "p_0"] \arrow[dd, "s_0", swap] && X_0 \arrow[dd, "f_0"]
            \\
            \\
            B_0 \arrow[uurr, "\exists ! u_0", dashed] \arrow[rr, "q_0", swap] && Y_0&{}
        \end{tikzcd}\begin{tikzcd}
        B_1 \arrow[rrd, "q_1", bend left] \arrow[d, "(d_0 {,} d_1)", swap] \arrow[dr, "\exists ! u_1", dashed] & & \\
        B_0 \times B_0 \arrow[dr, "u_0 \times u_0", swap] & X_1 \arrow[r, "f_1"] \arrow[d, "(d_0 {,} d_1)", swap] \arrow[dr, phantom, "\lrcorner", very near start]  & Y_1 \arrow[d, "(d_0 {,} d_1)"] \\
        & X_0 \times X_0 \arrow[r, "f_0 \times f_0", swap] & Y_0 \times Y_0
    \end{tikzcd}
$$
 By construction, $u : = (u_0, u_1): \mathbb{B} \to \X$ is a morphism of graphs. We show that $u : = (u_0, u_1): \B \to \X$ is a functor. Fix $k \in \{0, 1\}$. Then $u: \B \to \X$ respects identities by the universal property of $X_1$, as witnessed by the following commutative diagrams.

\begin{equation*}
    \begin{tikzcd}
        B_0 \arrow[r, "u_0"] \arrow[dd, "i", swap] \arrow[dr, "q_0"] & X_0 \arrow[r, "i"] \arrow[d, "f_0"] & X_1 \arrow[dd, "f_1"]  \\
        & X_0 \arrow[dr, "i"] & \\
        B_1 \arrow[r, "u_1", swap] \arrow[rr, "q_1", bend left] & X_1 \arrow[r, "f_1", swap] & Y_1&{}
    \end{tikzcd}
    \begin{tikzcd}
        B_0 \arrow[r, "u_0"] \arrow[dd, "i", swap] \arrow[dr, "1_{B_0}"] & X_0 \arrow[r, "i"] \arrow[ddr, "1_{X_0}"]& X_1 \arrow[dd, "d_k"]  \\
        & B_0 \arrow[dr, "u_0", swap] & \\
        B_1 \arrow[r, "u_1", swap] \arrow[ur, "d_k"]  & X_1 \arrow[r, "d_k", swap] & X_0
    \end{tikzcd}
\end{equation*}
      Similarly, the following commutative diagrams show that it respects composition.

    \begin{equation*}
         \begin{tikzcd}[row sep = 35, column sep = 35]
        B_2 \arrow[r, "u_2"] \arrow[dd, "m", swap] \arrow[dr, "q_2"] & X_2 \arrow[r, "m"] \arrow[d, "f_2"] & Y_2 \arrow[dd, "f_1"]  \\
        & X_0 \arrow[dr, "m"] & \\
        B_1 \arrow[r, "u_1", swap] \arrow[rr, "q_1", bend left] & X_1 \arrow[r, "f_1", swap] & Y_1&{}
    \end{tikzcd}
    \begin{tikzcd}
        B_2 \arrow[rr, "u_2"] \arrow[ddd, "m", swap] \arrow[dr, "\pi_k"] & & X_2 \arrow[r, "m"] \arrow[d, "\pi_k"] & X_1 \arrow[ddd, "d_k"]  \\
        & B_1 \arrow[r, "u_1"]  \arrow[d, "d_k"]& X_1 \arrow[ddr, "d_k"] & \\
        & B_0 \arrow[drr, "u_0"] & & \\
        B_1 \arrow[ur, "d_k"] \arrow[rr, "u_1", swap] & & X_1 \arrow[r, "d_k", swap] & X_0
    \end{tikzcd}
    \end{equation*}

    Hence $u: \B \to \X$ is an internal functor. But observe that $u_{0}: B_{0} \to X_{0}$ is the unique morphism satisfying $f_{0}u_{0} = q_{0}$ and $u_{0}s_{0} = p_{0}$, since $\Lagr \perp_{1} \R$ in $\E$. Moreover, it is clear by the construction of $u_{1}$ via the pullback that $(u_{0}, u_{1}): \B \to \X$ is the unique morphism of graphs providing a factorisation $fu=q$. But also $u_{1}s_{1}=p_{1}$, as per the following calculations using the universal property of $X_{1}$.

    $$\begin{tikzcd}
        A_{1}\arrow[dd, "p_{1}"']\arrow[rr, "s_{1}"] && B_{1}\arrow[rrdd, "q_{1}"']\arrow[rr, "u_{1}"] && X_{1}\arrow[dd, "f_{1}"]
        \\
        \\
        X_{1}\arrow[rrrr, "f_{1}"']&&&&Y_{1}&{}
    \end{tikzcd}
    \begin{tikzcd}
        A_1 \arrow[rr, "s_1"] \arrow[dr, "d_k"] \arrow[dd, "p_1"'] & &  B_1 \arrow[d, "d_k"] \arrow[r, "u_1"] & X_1 \arrow[dd, "d_k"] \\
        & A_0 \arrow[r, "s_0"] \arrow[drr, "p_0"'] & B_0 \arrow[dr, "u_0"] & \\
        X_1 \arrow[rrr, "d_k"'] & & & X_0
    \end{tikzcd}$$
   
    Thus $\Lagr' \perp_{1} \R'$. For the two-dimensional aspect of orthogonality, let $fp^0=q^0s$, $fp^1=q^1s$, $\overline{\alpha}: p^0  \Rightarrow p^1$ and $\overline{\beta}: q^0 \Rightarrow q^1$ be internal natural transformations satisfying $f.\overline{\alpha} = \overline{\beta}. s$, and let $u^0: \B \to \X$ and $u^1: \B \to \X$ be the uniquely induced maps from the one-dimensional aspect of orthogonality. Then by fully faithfulness of $\CatE(\B, f): \CatE(\B, \X) \to \CatE(\B, \Y)$, there is a unique internal natural transformation $\overline{\gamma}: u^{0} \Rightarrow u^{1}$ satisfying $f.\overline{\gamma} = \overline{\beta}$. As such, the components assigner for $\overline{\gamma}$ is induced by the universal property of $X_{1}$ as displayed below.

   $$\begin{tikzcd}
        B_0 \arrow[rrd, "\beta", bend left] \arrow[ddr, bend right = 30, "(u^0_0 {,} u^1_{0})", swap] \arrow[dr, "\exists ! \gamma", dashed] & & 
        \\
         & X_1 \arrow[r, "f_1"] \arrow[d, "(d_0 {,} d_1)", swap] \arrow[dr, phantom, "\lrcorner", very near start]  & Y_1 \arrow[d, "(d_0 {,} d_1)"] \\
        & X_0 \times X_0 \arrow[r, "f_0 \times f_0", swap] & Y_0 \times Y_0
    \end{tikzcd}$$

     Finally, the following diagrams for $k \in \{0, 1\}$ verify that $\overline{\gamma}.s = \overline{\alpha}$, completing the proof.

    \begin{equation*}
        \begin{tikzcd}
            A_0 \arrow[r, "s_0"] \arrow[d, "\alpha"'] \arrow[drr, "p^k_0"'] & B_0 \arrow[r, "\gamma"] \arrow[dr, "u^k_0"] & X_1 \arrow[d, "d_k"] \\
            X_1 \arrow[rr, "d_k"'] & & X_0
        \end{tikzcd}
        \qquad
        \begin{tikzcd}
            A_0 \arrow[r, "s_0"] \arrow[d, "\alpha"'] & B_0 \arrow[r, "\gamma"] \arrow[dr, "\beta"'] & X_1 \arrow[d, "f_1"] \\
            X_1 \arrow[rr, "f_1"'] & & Y_1
        \end{tikzcd}
    \end{equation*}
\end{proof}

\begin{lem}\label{lemma factorisation}
    Any internal functor $f: \X \to \Y$ may be factorised as $f=rl$ with $l \in \Lagr'$ and $r\in\R'$.
\end{lem}

\begin{proof}
     Let $f: \X \to \Y$ in $\CatE$. Using the $(\Lagr, \R)$ orthogonal factorisation system, we obtain a unique factorisation of $f_0$ in $\E$ depicted below left. We construct $C_1$  and maps $r_1 : C_1 \to Y_1$ and $(d_0, d_1): C_1 \to C_0 \times C_0$ via the pullback in $\E$ depicted below right.

    \begin{equation}\label{construction of r and l 0}
    \begin{tikzcd}
        X_0 \arrow[rr, "f_0"]\arrow[dr, "l_0", swap] & & Y_0&{} \\
        & C_0 \arrow[ur, "r_0"']
    \end{tikzcd} \begin{tikzcd}
        C_1 \arrow[r, "r_1"] \arrow[d, "(d_0{,} d_1)", swap] \arrow[dr, phantom, "\lrcorner", very near start] & Y_1 \arrow[d, "(d_0{,} d_1)"]   \\
        C_0 \times C_0 \arrow[r, "r_0 \times r_0", swap] & Y_0 \times Y_0
    \end{tikzcd}
    \end{equation}

 Define $l_1: X_1 \to C_1$ by the universal property of this pullback, as depicted below.

$$\begin{tikzcd}
        X_1 \arrow[rrd, "f_1", bend left] \arrow[d, "(d_0 {,} d_1)"'] \arrow[dr, "\exists !l_1", dashed] & & \\
        X_0 \times X_0 \arrow[dr, "l_0 \times l_0"'] & C_1 \arrow[r, "r_1"] \arrow[d, "(d_0 {,} d_1)"] \arrow[dr, phantom, "\lrcorner", very near start]  & Y_1 \arrow[d, "(d_0 {,} d_1)"] \\
        & C_0 \times C_0 \arrow[r, "r_0 \times r_0"'] & Y_0 \times Y_0
    \end{tikzcd}$$

 Then $(f_{0}, f_{1}) = (r_{0}, r_{1})\circ (l_{0}, l_{1})$ is clearly a factorisation at the level of morphisms of graphs. It remains to give an internal category structure to the graph in $\E$ displayed below, and to show that these morphisms of graphs are well-defined as internal functors. Once we have shown this, it will follow by construction that $l \in \Lagr'$ and $r \in \R'$.

$$\C := \begin{tikzcd}
    C_{1}\arrow[r, shift left = 2, "d_{0}"]\arrow[r, shift right =2, "d_{1}"'] &C_{0}
\end{tikzcd}$$

 Define the identity assigner $i: C_0 \to C_1$ for $\C$ using the universal property of $C_1$, as depicted below left. Then construct $C_{2} \in \E$ as the pullback depicted below right.

\begin{equation}\label{construction of i}
    \begin{tikzcd}
        C_0 \arrow[r, "r_0"] \arrow[dr, "i", dashed] \arrow[ddr, bend right, "(1_{C_{0}}{,}1_{C_{0}})", swap] & Y_0 \arrow[dr, "i"] & \\
        & C_1 \arrow[r, "r_1"] \arrow[d, "(d_0{,} d_1)"] \arrow[dr, phantom, "\lrcorner", very near start] & Y_1 \arrow[d, "(d_0{,} d_1)"]  \\
        & C_0 \times C_0 \arrow[r, "r_0 \times r_0", swap] & Y_0 \times Y_0 &{}
    \end{tikzcd}\begin{tikzcd}
        C_2 \arrow[rr, "\pi_0"] \arrow[dd, "\pi_1", swap] \arrow[dr, phantom, "\lrcorner", very near start] && C_1 \arrow[dd, "d_1"]
        \\&{}
        \\
        C_1 \arrow[rr, "d_0", swap] && C_0
    \end{tikzcd}
\end{equation}

 Define $r_2: C_2 \to Y_2$ by the universal property of $Y_2$, as described in Remark~\ref{Explicit internal functor}. Then define $m: C_2 \to C_1$ by the universal property of $C_1$ as depicted below left, given the commutativity of the diagram depicted below right.

\begin{equation}\label{construction of m}
    \begin{tikzcd}[column sep = small]
        C_2 \arrow[r, "r_2"] \arrow[dr, "m", dashed] \arrow[d, "(\pi_0 {,} \pi_1)", swap] & Y_2 \arrow[dr, "m"] & \\
       C_1 \times C_1 \arrow[dr, "d_0 \times d_1", swap] & C_1 \arrow[r, "r_1"] \arrow[d, "(d_0{,} d_1)"] \arrow[dr, phantom, "\lrcorner", very near start] & Y_1 \arrow[d, "(d_0{,} d_1)"]  \\
        & C_0 \times C_0 \arrow[r, "r_0 \times r_0", swap] & Y_0 \times Y_0&{}
    \end{tikzcd}\begin{tikzcd}[column sep = small]
        C_2 \arrow[r, "r_2"] \arrow[d, "(\pi_0 {,} \pi_1)", swap]
        & Y_2 \arrow[dr, "m"]\arrow[d, "(\pi_{0}{,}\pi_{1})"']  
        \\
       C_1 \times C_1 \arrow[dr, "d_0 \times d_1", swap] \arrow[r, "r_{1}\times r_{1}"]
       & Y_1 \times Y_{1}\arrow[dr, "(d_{0}{,}d_{1})"description]
       & Y_1 \arrow[d, "(d_0{,} d_1)"]
       \\
        & C_0 \times C_0 \arrow[r, "r_0 \times r_0", swap] 
        & Y_0 \times Y_0
    \end{tikzcd}
\end{equation}

 We now consider the internal category axioms for $\mathbb{C}$. Sources and targets for identities and composites hold by construction. This allows us to define the maps $\pi_{1, 3}, m_{1}, m_{0}, \pi_{0, 3}: C_{3} \to C_{2}$ and $i_{0}, i_{1}: C_{1} \to C_{2}$ as in Remark~\ref{internal category long definition}. Furthermore, define $r_3: C_3 \to Y_3$ in the obvious way, using the universal property of $Y_3$. It remains to check the associativity law and the left and right unit laws.

 To check associativity, we use the universal property of $C_1$, and the defining properties of $m$ and the relevant pullbacks. For $k \in \{0,1\}$ and $j=k+1$ mod $2$, we have: 
\begin{equation*}
    \begin{tikzcd}
       C_3 \arrow[rr, "m_j"] \arrow[d, "m_k"'] \arrow[dr, "\pi_{3,j}"] & & C_2 \arrow[d,"\pi_k"] \arrow[r, "m"] & C_1 \arrow[ddd, "d_k"] 
       \\
        C_2 \arrow[dr, "\pi_k"'] \arrow[dd, "m", swap] & C_2 \arrow[d, "\pi_k"] \arrow[r, "m"] & C_1 \arrow[ddr, "d_k"]& \\
        & C_1  \arrow[drr, "d_k"]& \\
       C_1 \arrow[rrr, "d_k"'] & &   & C_0 &{}
    \end{tikzcd}
     \begin{tikzcd}
        C_3
        \arrow[rr, "m_1"] \arrow[d, "m_0"'] \arrow[dr, "r_3"]
        & & C_2 \arrow[d,"r_2"] \arrow[r, "m"] 
        & C_1 \arrow[ddd, "r_1"] \\
        C_2  \arrow[dd, "m"'] \arrow[dr, "r_2"'] 
        & Y_3 \arrow[r, "m_1"] \arrow[d, "m_0"]
        & Y_2 \arrow[ddr, "m"] \\
        &  Y_2 \arrow[drr,"m"] 
        & \\
        C_1\arrow[rrr, "r_1"'] 
        & &
        & Y_1
    \end{tikzcd}
\end{equation*}

 We now consider the unit laws. We first note that the equations $r_{2}.i_{k}^{\mathbb{C}} = i_{k}^{\mathbb{Y}}.r_{1}$ for $k \in \{0, 1\}$ hold by the universal property of $Y_{2}$, as per the following calculations.
 
 \begin{equation*}
    \begin{tikzcd}
        C_1 \arrow[rr, "i_k"] \arrow[dr, "d_k"] \arrow[dd, "r_1"'] &  & C_2 \arrow[r, "r_2"]  \arrow[d, "\pi_k"] & Y_2 \arrow[ddd, "\pi_k"] \\
         & C_0 \arrow[d, "r_0"] \arrow[r, "i"] & C_1 \arrow[ddr, "r_1"]  & \\
         Y_1 \arrow[r, "d_k"'] \arrow[d, "i_k"'] & Y_0 \arrow[drr, "i"] & & \\
         Y_2 \arrow[rrr, "\pi_k"'] &  & & Y_1&{}
    \end{tikzcd}
    \begin{tikzcd}[row sep = 37, column sep = 37]
        C_1 \arrow[r, "i_k"] \arrow[d, "r_1"'] \arrow[dr, "1_{C_1}"]& C_2 \arrow[d, "\pi_j"]\arrow[r, "r_2"] & Y_2 \arrow[dd, "\pi_j"] \\
         Y_1  \arrow[d, "i_k"'] \arrow[drr, "1_{C_1}"'] & C_1 \arrow[dr, "r_1"]  & \\
         Y_2 \arrow[rr, "\pi_j"'] &   & Y_1
    \end{tikzcd}
\end{equation*}
 
 The left and right unit laws for $\C$ hence follow from the universal property of $C_{1}$, given the calculations displayed below where $j = k+1 \text{ mod } 2$.

\begin{equation*}
    \begin{tikzcd}
        C_1 \arrow[rr, "i_k"] \arrow[dd, "1_{C_1}"'] \arrow[dr, "d_k"'] &  & C_2 \arrow[r, "m"] \arrow[d, "\pi_k"] & C_1 \arrow[dd, "d_k"] \\
        & C_0 \arrow[r, "i"] \arrow[drr, "1_{C_0}"'] & C_1 \arrow[dr, "d_k"] & \\
        C_1 \arrow[rrr, "d_k"'] & &  & C_0&{}
    \end{tikzcd}\begin{tikzcd}
        C_1 \arrow[rr, "i_k"] \arrow[dd, "1_{C_1}"'] \arrow[dr, "r_1"'] &  & C_2 \arrow[r, "m"] \arrow[d, "r_2"] & C_1 \arrow[dd, "r_1"] \\
        & Y_1 \arrow[r, "i_k"] \arrow[drr, "1_{Y_1}"'] & Y_2 \arrow[dr, "m"] & \\
        C_1 \arrow[rrr, "r_1"'] & &  & Y_1&{}
    \end{tikzcd}
\end{equation*}
\begin{equation*}
    \begin{tikzcd}
        C_1 \arrow[r, "i_k"] \arrow[dd, "1_{C_1}"'] \arrow[dr, "1_{C_1}"'] & C_2 \arrow[d, "\pi_j"] \arrow[r, "m"] & C_1 \arrow[dd, "d_j"] \\
        & C_1 \arrow[dr, "d_j"] & \\
        C_1 \arrow[rr, "d_j"'] & & C_0
    \end{tikzcd}
\end{equation*}

 So $\C$ is a category internal to $\E$. It is clear from the construction of the identity assigner in \autoref{construction of i} and composition in \autoref{construction of m} that the morphism of graphs $r := (r_0, r_1) : \C \to \Y$ is well-defined as an internal functor, which is moreover evidently fully faithful and $\R$-on-objects as per \autoref{construction of r and l 0}.

 It remains to show that the morphism of graphs $l := (l_0, l_1): \X \to \C$ is well-defined as an internal functor. Once again, we do this using the universal property of $C_{1}$ as a pullback. Define $l_2: X_2 \to C_2$  by the universal property of $C_2$, as described in Remark~\ref{Explicit internal functor}. Fix $k \in \{0, 1\}$ as above. Respect for identities for $l: \X \to \C$ is exhibited by the commutativity of the diagrams in $\E$ displayed below.

\begin{equation*}
    \begin{tikzcd}
        X_0 \arrow[r, "l_0"] \arrow[d, "i"'] \arrow[dr, "1_{X_0}"] & C_0 \arrow[r, "i"] \arrow[ddr, "1_{C_0}"] & C_1 \arrow[dd, "d_k"] \\
        X_1 \arrow[d, "l_1"'] \arrow[r, "d_k"] & X_0 \arrow[dr, "l_0"'] & \\
        C_1 \arrow[rr, "d_k"'] & & C_0
    \end{tikzcd}
    \qquad
    \begin{tikzcd}
        X_0 \arrow[r, "l_0"] \arrow[d, "i"'] \arrow[dr, "f_0"] & C_0 \arrow[r, "i"] \arrow[d, "r_0"] & C_1 \arrow[dd, "r_1"] \\
        X_1 \arrow[d, "l_1"'] \arrow[drr, "f_1"'] & Y_0 \arrow[dr, "i"] & \\
        C_1 \arrow[rr, "r_1"'] & & Y_1
    \end{tikzcd}
\end{equation*}
 Finally, respect for composition for $l: \X \to \C$ follows from the commutativity of the diagrams in $\E$ displayed below. This completes the proof.

\begin{equation*}
    \begin{tikzcd}
        X_2 \arrow[rr, "l_2"] \arrow[dr, "\pi_k"] \arrow[dd, "m"'] & & C_2 \arrow[d, "\pi_k"] \arrow[r, "m"] & C_1 \arrow[ddd, "d_k"] \\
        & X_1 \arrow[r, "l_1"'] \arrow[d, "d_k"] & C_1 \arrow[ddr, "d_k"] & \\
        X_1 \arrow[r, "d_k"] \arrow[d, "l_1"'] & X_0 \arrow[drr, "l_0"] & & \\
        C_1 \arrow[rrr, "d_k"'] & & & C_0
    \end{tikzcd}
    \qquad
    \begin{tikzcd}[row sep = 38, column sep = 38]
        X_2 \arrow[r, "l_2"] \arrow[d, "m"'] \arrow[dr, "f_2"] & C_2 \arrow[d, "r_2"] \arrow[r, "m"] & C_1 \arrow[dd, "r_1"] \\
        X_1 \arrow[drr, "f_1"'] \arrow[d, "l_1"']& Y_2 \arrow[dr, "m"] & \\
        C_1 \arrow[rr, "r_1"'] &  & Y_1 
    \end{tikzcd}
\end{equation*}

\end{proof}


\begin{thebibliography}{EKVdL05}

\bibitem[BE72]{bastiani1972categories}
Andr{\'e}e Bastiani and Charles Ehresmann.
\newblock Categories of sketched structures.
\newblock {\em Cahiers de topologie et geometrie differentielle}, 13(2):104--214, 1972.

\bibitem[BG14]{bourke2014two}
John Bourke and Richard Garner.
\newblock Two-dimensional regularity and exactness.
\newblock {\em Journal of Pure and Applied Algebra}, 218(7):1346--1371, 2014.

\bibitem[B{\"o}h20]{bohm2020gray}
Gabriella B{\"o}hm.
\newblock The gray monoidal product of double categories.
\newblock {\em Applied Categorical Structures}, 28(3):477--515, 2020.

\bibitem[Bor94]{borceux1994handbook}
Francis Borceux.
\newblock {\em Handbook of categorical algebra: volume 1, Basic category theory}, volume~1.
\newblock Cambridge University Press, 1994.

\bibitem[Bou77]{bousfield1977constructions}
Aldridge~K Bousfield.
\newblock Constructions of factorization systems in categories.
\newblock {\em Journal of Pure and Applied Algebra}, 9(2-3):207--220, 1977.

\bibitem[Bou10]{bourke2010codescent}
John Bourke.
\newblock {\em Codescent objects in 2-dimensional universal algebra}.
\newblock PhD thesis, University of Sydney, 2010.

\bibitem[CLW93]{CARBONI1993145}
Aurelio Carboni, Stephen Lack, and R.F.C. Walters.
\newblock Introduction to extensive and distributive categories.
\newblock {\em Journal of Pure and Applied Algebra}, 84(2):145--158, 1993.

\bibitem[Day72]{day1972reflection}
Brian Day.
\newblock A reflection theorem for closed categories.
\newblock {\em Journal of pure and applied algebra}, 2(1):1--11, 1972.

\bibitem[Day06]{day2006adjoint}
Brian Day.
\newblock On adjoint-functor factorisation.
\newblock In {\em Category Seminar: Proceedings Sydney Category Theory Seminar 1972/1973}, pages 1--19. Springer, 2006.

\bibitem[Dia75]{diaconescu1975axiom}
Radu Diaconescu.
\newblock Axiom of choice and complementation.
\newblock {\em Proceedings of the American Mathematical Society}, 51(1):176--178, 1975.

\bibitem[Ehr59]{ehresmann1959categories}
Charles Ehresmann.
\newblock {\em Cat{\'e}gories topologiques et cat{\'e}gories diff{\'e}rentiables}.
\newblock Librairie universitaire, 1959.

\bibitem[Ehr63]{ehresmann1963categories}
Charles Ehresmann.
\newblock Cat{\'e}gories structur{\'e}es.
\newblock In {\em Annales scientifiques de l'{\'E}cole Normale Sup{\'e}rieure}, volume~80, pages 349--426, 1963.

\bibitem[EKVdL05]{everaert2005model}
Tomas Everaert, R.W.~Kieboom, and Tim Van~der Linden.
\newblock Model structures for homotopy of internal categories.
\newblock {\em Theory Appl. Categ}, 15(3):66--94, 2005.

\bibitem[Fre19]{frey2019characterizing}
Jonas Frey.
\newblock Characterizing partitioned assemblies and realizability toposes.
\newblock {\em Journal of Pure and Applied Algebra}, 223(5):2000--2014, 2019.

\bibitem[FS90]{freyd1990categories}
Peter~J. Freyd and Andre Scedrov.
\newblock {\em Categories, allegories}.
\newblock Elsevier, 1990.

\bibitem[GJ09]{goerss2009simplicial}
Paul~G. Goerss and John~F. Jardine.
\newblock {\em Simplicial homotopy theory}.
\newblock Springer Science \& Business Media, 2009.

\bibitem[Gro60]{grothendieck1960techniques}
Alexander Grothendieck.
\newblock Techniques de construction et th{\'e}oremes d’existence en g{\'e}om{\'e}trie alg{\'e}brique. iv. les sch{\'e}mas de hilbert.
\newblock {\em S{\'e}minaire Bourbaki}, 6(221):249--276, 1960.

\bibitem[HM]{HughesMiranda20242categories}
Calum Hughes and Adrian Miranda.
\newblock $2$-categories of categories, discrete opfibration classifiers and the axiom of replacement.
\newblock In preparation.

\bibitem[JM95]{joyal1995algebraic}
Andr{\'e} Joyal and Ieke Moerdijk.
\newblock {\em Algebraic set theory}, volume 220.
\newblock Cambridge University Press, 1995.

\bibitem[Joh02]{johnstone2002sketches}
Peter~T. Johnstone.
\newblock {\em Sketches of an Elephant: A Topos Theory Compendium, Volume 1}.
\newblock Oxford University Press, 2002.

\bibitem[Joh14]{johnstone2014topos}
Peter~T. Johnstone.
\newblock {\em Topos theory}.
\newblock Courier Corporation, 2014.

\bibitem[JT06]{joyal2006strong}
Andr{\'e} Joyal and Myles Tierney.
\newblock Strong stacks and classifying spaces.
\newblock In {\em Category Theory: Proceedings of the International Conference held in Como, Italy, July 22--28, 1990}, pages 213--236. Springer, 2006.

\bibitem[JW78]{johnstonealgebraic}
Peter~T. Johnstone and Gavin~C Wraith.
\newblock Algebraic theories in toposes.
\newblock In {\em Indexed Categories and Their Applications}. Springer, 1978.

\bibitem[JY21]{johnson20212}
Niles Johnson and Donald Yau.
\newblock {\em 2-dimensional categories}.
\newblock Oxford University Press, USA, 2021.

\bibitem[Koc06]{kock2006synthetic}
Anders Kock.
\newblock {\em Synthetic differential geometry}, volume 333.
\newblock Cambridge University Press, 2006.

\bibitem[Lac09]{lack20092}
Stephen Lack.
\newblock A 2-categories companion.
\newblock In {\em Towards higher categories}, pages 105--191. Springer, 2009.

\bibitem[Lan17]{landry2017categories}
Elaine Landry.
\newblock {\em Categories for the working philosopher}.
\newblock Oxford University Press, 2017.

\bibitem[Law]{Lawvere2016Grothendieck}
F.~William Lawvere.
\newblock Alexander {G}rothendieck and the concept of space.
\newblock Address, CT15 Aveiro 2016.

\bibitem[Law63]{lawvere1963functorial}
F.~William Lawvere.
\newblock Functorial semantics of algebraic theories.
\newblock {\em Proceedings of the National Academy of Sciences}, 50(5):869--872, 1963.

\bibitem[Law64]{lawvere1964elementary}
F.~William Lawvere.
\newblock An elementary theory of the category of sets.
\newblock {\em Proceedings of the national academy of sciences}, 52(6):1506--1511, 1964.

\bibitem[Law66]{lawvere1966category}
F.~William Lawvere.
\newblock The category of categories as a foundation for mathematics.
\newblock In {\em Proceedings of the Conference on Categorical Algebra: La Jolla 1965}, pages 1--20. Springer, 1966.

\bibitem[Lei14]{leinster2014rethinking}
Tom Leinster.
\newblock Rethinking set theory.
\newblock {\em The American Mathematical Monthly}, 121(5):403--415, 2014.

\bibitem[LM05]{lawvere2005elementary}
F.~William Lawvere and Colin McLarty.
\newblock An elementary theory of the category of sets (long version) with commentary.
\newblock {\em Reprints in Theory and Applications of Categories}, 11:1--35, 2005.

\bibitem[Lur09]{lurie2009higher}
Jacob Lurie.
\newblock {\em Higher topos theory}.
\newblock Princeton University Press, 2009.

\bibitem[Mak96]{makkai1996avoiding}
Michael Makkai.
\newblock Avoiding the axiom of choice in general category theory.
\newblock {\em Journal of pure and applied algebra}, 108(2):109--173, 1996.

\bibitem[Mir18]{miranda2022internal}
Adrian Miranda.
\newblock Internal categories.
\newblock Master's thesis, Macquarie University, 2018.
\newblock Available at \url{https://figshare.mq.edu.au/articles/thesis/Internal_categories/19434626/1}.

\bibitem[ML13]{mac2013categories}
Saunders Mac~Lane.
\newblock {\em Categories for the working mathematician}, volume~5.
\newblock Springer Science \& Business Media, 2013.

\bibitem[MM12]{maclane2012sheaves}
Saunders MacLane and Ieke Moerdijk.
\newblock {\em Sheaves in geometry and logic: A first introduction to topos theory}.
\newblock Springer Science \& Business Media, 2012.

\bibitem[{nLa}23]{nlab:fully_formal_etcs}
{nLab authors}.
\newblock fully formal {{E}}{{T}}{{C}}{{S}}.
\newblock \url{https://ncatlab.org/nlab/show/fully+formal+ETCS}, December 2023.
\newblock \href{https://ncatlab.org/nlab/revision/fully+formal+ETCS/31}{Revision 31}.

\bibitem[Osi74]{osius1974categorical}
Gerhard Osius.
\newblock Categorical set theory: a characterization of the category of sets.
\newblock {\em Journal of Pure and Applied Algebra}, 4(1):79--119, 1974.

\bibitem[PR91]{power1991characterization}
John Power and Edmund Robinson.
\newblock A characterization of pie limits.
\newblock In {\em Mathematical Proceedings of the Cambridge Philosophical Society}, volume 110, pages 33--47. Cambridge University Press, 1991.

\bibitem[Rob12]{roberts2012internal}
David~Michael Roberts.
\newblock Ineternal categories, anafunctors and localisations.
\newblock {\em Theory and Applications of Categories}, 26(29):788--829, 2012.

\bibitem[Rob21]{roberts2021elementary}
David~Michael Roberts.
\newblock The elementary construction of formal anafunctors.
\newblock {\em Categories and General Algebraic Structures with Applications}, 15(1):183--229, 2021.

\bibitem[RV22]{riehl2022elements}
Emily Riehl and Dominic Verity.
\newblock {\em Elements of $\infty$-category theory}, volume 194.
\newblock Cambridge University Press, 2022.

\bibitem[Ste24]{stenzel2024infty2category}
Raffael Stenzel.
\newblock The $(\infty,2)$-category of internal $(\infty,1)$-categories, arXiv preprint arXiv:2402.01396, 2024.

\bibitem[Str76]{street1976limits}
Ross Street.
\newblock Limits indexed by category-valued 2-functors.
\newblock {\em Journal of Pure and Applied Algebra}, 8(2):149--181, 1976.

\bibitem[Str80]{street1980cosmoi}
Ross Street.
\newblock Cosmoi of internal categories.
\newblock {\em Transactions of the American Mathematical Society}, 258(2):271--318, 1980.

\bibitem[Str82]{street1982two}
Ross Street.
\newblock Two-dimensional sheaf theory.
\newblock {\em Journal of Pure and Applied Algebra}, 23(3):251--270, 1982.

\bibitem[Str06]{street2006elementary}
Ross Street.
\newblock Elementary cosmoi i.
\newblock In {\em Category Seminar: Proceedings Sydney Category Theory Seminar 1972/1973}, pages 134--180. Springer, 2006.

\bibitem[Web07]{weber2007yoneda}
Mark Weber.
\newblock Yoneda structures from 2-toposes.
\newblock {\em Applied Categorical Structures}, 15:259--323, 2007.

\end{thebibliography}
\end{document}